\magnification=1150
\vsize=185mm
\hsize=135mm
\voffset=5mm
\null


\def\equp{1.1}
\def\Qbddu{1.2}

\def\SdecayHaaa{2.1}
\def\SdecayHba{2.2}
\def\Sdecayu{2.3}
\def\hypdecayuzeto{2.4}
\def\SdecayHaa{2.5}
\def\SdecayHb{2.6}
\def\DefLambdaStar{2.7}
\def\Gopen{2.8}
\def\GopenA{2.9}
\def\GopenB{2.10}
\def\GopenC{2.11}

\def\LinearBehavior{3.1}
\def\LinearBehaviorB{3.2}
\def\AprioriEst{3.3}

\def\DefMorrey{4.1}
\def\PropMorreyA{4.2}
\def\PropMorreyB{4.3}
\def\RegulMorrey{4.4}
\def\ContractionMorrey{4.5}
\def\WeightedLinfty{4.6}
\def\compMorreyGnormA{4.7}
\def\compMorreyGnormB{4.8}
\def\compMorreyGnorm{4.9}
\def\compMorreyintermed{4.10}
\def\SGKdefSelfsimilA{4.11}
\def\SGKdefSelfsimilB{4.12}
\def\SGKdefrescaledOmega{4.13}
\def\QGKeqrho{4.14}
\def\SGKdefenergyOmega{4.15}
\def\QGKidA{4.16}
\def\QGKidB{4.17}
\def\QGKenergypos{4.18}
\def\QGKLtwocontrol{4.19}

\def\RHSa{5.1}
\def\RHSb{5.2}

\def\DefScaling{6.1}
\def\SupMorrey{6.2}
\def\SMorreyInvEstA{6.3}
\def\SMorreyInvEstAA{6.4}
\def\SMorreyInvEstAB{6.5}
\def\intinv{6.6}
\def\intinvB{6.7}
\def\SMorreyInvEst{6.8}
\def\SMorreyInvEstB{6.9}
\def\Varconsthalf{6.10}

\def\Boundtzero{7.1}
\def\TEnegA{7.2}
\def\TEnegB{7.3}
\def\TEneg{7.4}

\def\depcontMorrey{8.1}
\def\supsbetaA{8.2}
\def\supsbeta{8.3}
\def\upvprp{8.4}
\def\supsbetaB{8.5}
\def\supsbetaBB{8.6}
\def\supsbetaC{8.7}

\def\ClaimKatoMorrey{A.1}
\def\SGKdds{A.2}
\def\SGKconvex{A.3}
\def\vjrho{A.4}
\def\vjSobolev{A.5}
\def\DecayGrad{A.6}
\def\varconstNabla{A.7}
\def\varconstNablaB{A.8}

\def\BKZ{1}
\def\BS{2}
\def\CFG{3}
\def\CDZ{4}
\def\Fi{5}
\def\FK{6}
\def\FM{7}
\def\GK{8}
\def\GV{9}
\def\Gig{10}
\def\GKa{11}
\def\GKb{12}
\def\GMi{13}
\def\HW{14}
\def\Kar{15}
\def\Kato{16}
\def\Kav{17}
\def\Kaw{18}
\def\KS{19}
\def\LN{20}
\def\LR{21}
\def\MMCPAM{22}
\def\MMJFA{23}
\def\MMJFAb{24}
\def\Miz{25}
\def\MizGlobal{26}
\def\NST{27}
\def\OS{28}
\def\Pac{29}
\def\PQS{30}
\def\PY{31}
\def\PYb{32}
\def\Qui{33}
\def\QuiB{34}
\def\QSb{35}
\def\SnTW{36}
\def\SouCRAS{37}
\def\SouCPDE{38}
\def\SouJFA{39}
\def\Tay{40}
\def\Tri{41}
\def\YZ{42}
\def\Wei{43}

\def\eps{\varepsilon}

\null

\newtoks\hautpagegauche
\newtoks\hautpagedroite
\newtoks\paragraphecourant
\newtoks\chapitrecourant
\hautpagegauche={}
\hautpagedroite={}
\headline={\ifodd\pageno\the\hautpagedroite\else\the\hautpagegauche\fi}

\font\TenEns=msbm10
\font\SevenEns=msbm7
\font\FiveEns=msbm5
\newfam\Ensfam
\def\Ens{\fam\Ensfam\TenEns}
\textfont\Ensfam=\TenEns
\scriptfont\Ensfam=\SevenEns
\scriptscriptfont\Ensfam=\FiveEns
\def\R{{\Ens R}}

\def\Rn{{\R}^n}

\def\eps{\varepsilon}


\font\itsmall=cmsl9

\font\eightrm=cmr9
\font\sixrm=cmr6
\font\fiverm=cmr5

\font\eighti=cmmi9
\font\sixi=cmmi6
\font\fivei=cmmi5

\font\eightsy=cmsy9
\font\sixsy=cmsy6
\font\fivesy=cmsy5

\font\eightit=cmti9
\font\eightsl=cmsl9
\font\eighttt=cmtt9

\def\eightpoint{\def\rm{\fam0\eightrm}
\textfont0=\eightrm
\scriptfont0=\sixrm
\scriptscriptfont0=\fiverm

\textfont1=\eighti
\scriptfont1=\sixi
\scriptscriptfont1=\fivei

\textfont2=\eightsy
\scriptfont2=\sixsy
\scriptscriptfont2=\fivesy

\textfont3=\tenex
\scriptfont3=\tenex
\scriptscriptfont3=\tenex

\textfont\itfam=\eightit \def\it{\fam\itfam\eightit}
\textfont\slfam=\eightsl \def\it{\fam\slfam\eightsl}
\textfont\ttfam=\eighttt \def\it{\fam\ttfam\eighttt}
}


\font\pc=cmcsc9
\font\itsmall=cmsl9

\def \trait (#1) (#2) (#3){\vrule width #1pt height #2pt depth #3pt}
\def \fin{\hfill
       \trait (0.1) (5) (0)
       \trait (5) (0.1) (0)
       \kern-5pt
       \trait (5) (5) (-4.9)
       \trait (0.1) (5) (0)
\medskip}
%

\paragraphecourant={\rmt MORREY SPACES AND CLASSIFICATION OF GLOBAL SOLUTIONS}
\chapitrecourant={\rmt Ph. SOUPLET}
\footline={\ifnum\folio=1 \hfill\folio\hfill\fi}
\hautpagegauche={\tenrm\folio\hfill\the\chapitrecourant\hfill}
\hautpagedroite={\ifnum\folio=1 \hfill\else\hfill\the\paragraphecourant\hfill\tenrm\folio\fi}

\font\rmb=cmbx8 scaled 1125 \rm

\centerline{\rmb MORREY SPACES AND CLASSIFICATION OF GLOBAL SOLUTIONS}
\vskip 0.5mm
\centerline{\rmb FOR A SUPERCRITICAL SEMILINEAR HEAT EQUATION IN $\R^n$}

\vskip 5mm
\centerline{\pc Philippe SOUPLET}
\vskip 3mm
\centerline{\itsmall Universit\'e Paris 13, Sorbonne Paris Cit\'e, CNRS UMR 7539}
\centerline{\itsmall Laboratoire Analyse G\'eom\'etrie et Applications}
\centerline{\itsmall 93430 Villetaneuse, France. Email: souplet@math.univ-paris13.fr}

\vskip 4mm

\baselineskip=12pt
\font\rmt=cmr9

\setbox1=\vbox{
\hsize=120mm
{\baselineskip=11pt \parindent=3mm \eightpoint \rmt
{\pc Abstract:}\ We prove the boundedness of global classical solutions for the semilinear heat equation
$u_t-\Delta u= |u|^{p-1}u$ in the whole space $\R^n$, with $n\ge 3$ and supercritical power $p>(n+2)/(n-2)$.
This is proved {\rmb without any radial symmetry or sign assumptions}, unlike in all the previously known results
for the Cauchy problem,
and under spatial decay assumptions on the initial data that are essentially optimal
in view of the known counter-examples. Moreover, we show that any global classical solution has to decay in time 
faster than $t^{-1/(p-1)}$, which is also optimal and in contrast with the subcritical case.

The proof relies on nontrivial modifications of techniques developed by 
Chou, Du and Zheng [\CDZ] and by Blatt and Struwe [\BS] for the case of convex bounded domains.
They are based on weighted energy estimates of Giga-Kohn type, combined
with an analysis of the equation in a suitable Morrey space. 
We in particular simplify the approach of [\BS] by establishing and using a result on global existence and decay 
for small initial data in the critical Morrey space $M^{2,4/(p-1)}(\Rn)$,
rather than $\eps$-regularity in a parabolic Morrey space. 
This method actually works for any convex, bounded or unbounded, smooth domain,
but at the same time captures some of the specific behaviors associated with the case of the whole space $\R^n$.

As a consequence we also prove that the set of initial data producing global solutions
 is open in suitable topologies,
and we show that the so-called ``borderline'' global weak solutions blow up in finite time 
and then become classical again and decay as $t\to\infty$.
All these results put into light the key role played by the Morrey space  $M^{2,4/(p-1)}$
in the understanding of the structure of the set of global solutions for $p>p_S$.

\vskip 0.2cm

{\pc Keywords:}\ Semilinear heat equation, supercritical power, boundedness of global solutions, decay,
Morrey spaces, weighted energy

}}

\hskip 2mm \hbox{\box1}

\bigskip

{\bf 1. Introduction.}
\medskip

In this article, we consider the following semilinear heat equation 
$$\left\{\eqalign{
u_t-\Delta u &= |u|^{p-1}u,\qquad x\in\Omega,\ t>0,\cr
\noalign{\vskip 0.5mm}
      u  &= 0,  \qquad x\in\partial\Omega,\ t>0,\cr
\noalign{\vskip 0.5mm}
       u(x,0) &= u_0(x),\qquad x\in\Omega.
}\right.
\eqno{(\equp)}
$$
Here $p>1$, $\Omega\subset \R^n$ is a possibly unbounded domain and $u_0\in L^\infty(\Omega)$.
Unless otherwise specified, it is assumed throughout that $\Omega$ is of class $C^{2+\alpha}$.
Problem (\equp) admits a unique, maximal classical solution.
Moreover, denoting its maximal existence time by $T_{max}=T_{max}(u_0)\in (0,\infty]$,
we have $\lim_{t\to T_{max}}\|u(t)\|_\infty=\infty$ whenever $T_{max}<\infty$.
For $1\le q\le \infty$, the norm in $L^q(\Omega)$
will be denoted by $\|\cdot\|_q$ 
(the domain will be omitted when no confusion arises).
 \medskip
 
Problem (\equp) is one of the best studied model cases in the theory of semilinear parabolic equations and it has received considerable 
attention in the last fifty years (see e.g.~the monograph [\QSb] for extensive references).
 In particular, it is well known that (\equp) admits both small data global solutions,
and finite time blowup solutions for large initial data $u_0$.
It is therefore natural to investigate whether or not it admits other kinds of solutions
(namely global unbounded classical solutions).
 
 \smallskip
 
Denote by $p_S$ the Sobolev critical exponent
$$p_S=(n+2)/(n-2)_+.$$
In the subcritical case $p<p_S$ with $\Omega$ bounded, it is known [\Fi] that any global solution 
is uniformly bounded, i.e.
$$
\sup_{t\geq 0}\|u(t)\|_\infty <\infty,
\eqno{(\Qbddu)}
$$
see~also [\Gig, \Kav, \SouCRAS, \Qui, \PQS, \QSb] and the references therein for further results such as a priori estimates or universal bounds and for 
similar properties for nonnegative global solutions of the Cauchy problem.
The critical and supercritical cases $p=p_S$ and $p>p_S$ (with $n\ge 3$) are more involved and not as well understood.
For a long time, the only available results concerned radially symmetric solutions
and were based mainly on intersection-comparison arguments (see~[\GV, \Miz, \MMCPAM, \MizGlobal, \CFG, \MMJFA]).
More recently, by introducing some new ideas, 
the following result about boundedness of global solutions for $p>p_S$ in a nonradial framework was obtained in [\CDZ]  and  [\BS]
(cf.~[\BS,~Proposition~6.6] for the general case; the special case $u_0\ge 0$ follows from the proofs of Proposition~3 and Theorem B in [\CDZ]).

\proclaim Theorem~A.
Assume $p>p_S$, $\Omega$ convex bounded and $u_0\in L^\infty(\Omega)$. 
If the solution $u$ of (\equp) is global, then property (\Qbddu) is true.
Moreover, 
$$\lim_{t\to\infty} \|u(t)\|_\infty=0.$$

\smallskip
The main goals of the present article are:
\smallskip
(a) to investigate the boundedness and decay of global classical solutions for the Cauchy problem
($\Omega=\R^n$), or more generally in unbounded domains.
We will consider also the related question of
describing the structure of the set of global solutions, including the so-called ``borderline'' weak solutions;
\smallskip
(b) to provide an alternative approach, motivated by those in [\CDZ] and [\BS], 
and in turn give a simpler proof of Theorem~A
(for bounded domains; as for the case of $\Rn$ or unbounded domains, it involves 
a wider variety of phenomena, and specific difficulties).

\medskip\medskip

{\bf 2. Main results.}
\medskip
\ \ {\bf 2.1. Boundedness and decay of global classical solutions.}
\medskip

Our main results are the following.

\proclaim Theorem~1.
Let $p>p_S$, $\Omega=\Rn$ and $u_0\in L^\infty(\Rn)$.
Assume that either
$$\nabla u_0\in L^q(\Rn)\quad\hbox{for some } q\in \Bigl[2,{n(p-1)\over p+1}\Bigr);
\eqno{(\SdecayHaaa)}$$
or
$$u_0\in BC^1(\Rn)\quad\hbox{and}\quad
|\nabla u_0(x)|=o\bigl(|x|^{-{2\over p-1}-1}\bigr)\quad \hbox{ as $|x|\to\infty$.}
\eqno{(\SdecayHba)}$$
If the solution $u$ of (\equp) is global, then property (\Qbddu) is true.
Moreover $u$ satisfies the decay property 
$$\lim_{t\to\infty} t^{1/(p-1)}\|u(t)\|_\infty=0.
\eqno{(\Sdecayu)}$$

\medskip

We point out that Theorem~1 extends to the nonradial framework a result of  [\MMCPAM] for radial solutions in $\Rn$
(see [\MMJFA, Theorem~5.14] and its proof).
\smallskip

We next consider general convex, bounded or unbounded, domains, for which we have the following results.
In all this article, we use the following:
\medskip

{\bf Notation.} The Gaussian heat kernel is denoted by 
$$G_t(x)=(4\pi t)^{-n/2}e^{-|x|^2/4t},\quad x\in \R^n,\ t>0.$$
Also, for any $\phi\in L^\infty(\Omega)$, we write $G_t*\phi:=G_t*\tilde\phi$,
where $\tilde\phi\in L^\infty(\R^n)$ is the extension of $\phi$ by $0$.
Also $BC^1(\Omega)$ denotes the set of functions of class $C^1$ in 
$\overline\Omega$ which are bounded along with their first order derivatives.

\bigskip

\proclaim Theorem~2.
Assume $p>p_S$, $\Omega$ convex and $u_0\in L^\infty(\Omega)$. 
If $\Omega$ is unbounded, assume in addition that $u_0\in BC^1(\Omega)$ and that
$$\lim_{t\to\infty} \Bigl(t^{p+1\over p-1}\|G_t*|\nabla u_0|^2\|_\infty+t^{2\over p-1}\|G_t*|u_0|^2\|_\infty\Bigr)=0.
\eqno{(\hypdecayuzeto)}$$
If the solution $u$ of (\equp) is global, then property (\Qbddu) is true.
Moreover $u$ satisfies the decay property (\Sdecayu).

\medskip

As we shall see, condition (\hypdecayuzeto) is for instance true under the more familiar assumptions:
$$|u_0|^{p+1}+|\nabla u_0|^2\in L^m(\Omega)\quad\hbox{for some } m\in \Bigl[1,{n\over 2}{p-1\over p+1}\Bigr),
\eqno{(\SdecayHaa)}$$
or $u_0\in BC^1(\Omega)$ and
$$|u_0(x)|+|x|{\hskip 1pt}|\nabla u_0(x)|=o\bigl(|x|^{-{2\over p-1}}\bigr)\quad \hbox{ as $|x|\to\infty$.}
\eqno{(\SdecayHb)}$$
We thus have:

\medskip

\proclaim Corollary~3.
Assume $p>p_S$, $\Omega$ convex and $u_0\in L^\infty(\Omega)$. If $\Omega$ is unbounded, assume in addition that
$u_0\in BC^1(\Omega)$ and that either (\SdecayHaa) or (\SdecayHb) is satisfied.
 If the solution $u$ of (\equp) is global, then property (\Qbddu) is true.
Moreover $u$ satisfies the decay property~(\Sdecayu).

\medskip

We stress that the above decay assumptions on $u_0$ are {\it not} technical. In fact, in the case $\Omega=\Rn$, the assumption (\SdecayHb) on the spatial decay of $u_0$, as well as the conclusion~(\Sdecayu) on the temporal decay of $u$, are essentially optimal; see~Remark~3.1.

\medskip\medskip

\goodbreak
\ \ {\bf 2.2. Borderline solutions and structure of the set of global solutions.}
\medskip

A related question is that of ``threshold'' or ``borderline'' solutions, introduced in~[\NST].
For $\Omega$ bounded and a given nontrivial $\phi\in L^\infty(\Omega)$, consider problem (\equp) with initial data $u_0=\lambda\phi$
and let 
$$\lambda^*=\sup\{\lambda>0;\, T_{max}(\lambda\phi)=\infty\}.
\eqno{(\DefLambdaStar)}$$
It is well known that if either $\phi\ge 0$ or $\phi\in H^1_0(\Omega)$, then $\lambda^*\in(0,\infty)$.
The behavior of the solution $u^*(t)=u(t;\lambda^*\phi)$ has been studied in many papers. One of the following three possibilities must occur:

\smallskip
(a) $u^*$ is global bounded;
\vskip 0.5mm
(b) $u^*$ is global unbounded;
\vskip 0.5mm
(c) $u^*$ blows up in finite time.
\smallskip

It is known [\Gig, \Qui] that (a) occurs if $p<p_S$. If $p=p_S$, $\Omega$ is a ball and $\phi$ is radial nonnegative, then (b) occurs [\GV].
If $p>p_S$ and $\Omega$ is convex then (c) occurs. This is a consequence of the fact that the set
$$\hbox{${\cal G}=\{u_0\in L^\infty(\Omega);\ T_{max}(u_0)=\infty\}$
is open in the $L^\infty$ topoplogy},
\eqno(\Gopen)$$
owing to [\BS]. 
Indeed, for any $u_0\in {\cal G}$,  $\|u(t)\|_\infty$ decays as $t\to\infty$ by [\BS, Proposition~6.6].  
Property (\Gopen) then follows easily 
from the continuous dependence in $L^\infty$ and the well-known fact that
the trivial solution is stable in $L^\infty$ (with $\Omega$ bounded).
Some results concerning the threshold behavior for $\Omega=\R^n$ and $p>p_S$ in the radial case 
can be found in [\Miz, \MizGlobal, \MMJFA, \MMJFAb, \QuiB]
(see also [\QSb]).

We here give the following extension of property (\Gopen) to the Cauchy problem (or to convex unbounded domains)
 in suitable topologies
(note that the $L^\infty$-topology is no longer relevant, see Remark~3.1(i)).
Namely we show that property (\Gopen) is true with respect to the critical Morrey norm $M^{2,4/(p-1)}$, 
as well as to the critival $L^q$ norm and to the
$L^\infty$ norm weighted by $|x|^{2/(p-1)}$, corresponding to the borderline decay.
All these topologies are quite ``natural'' in view of the scaling properties of the equation (see~Section~6)
and the result seems actually new even for convex bounded domains.

\proclaim Theorem 4.
Assume $p>p_S$, $\Omega$ convex and $u_0\in {\cal G}$.
\smallskip
(i) If $\Omega$ is unbounded, assume in addition that $u_0\in BC^1(\Omega)$ and that (\hypdecayuzeto) holds.
Then there exists $\eta=\eta(u_0)>0$ such that 
$$\bigl\{ v_0\in L^\infty(\Omega);\, 
\|u_0-v_0\|_{M^{2,4/(p-1)}} <\eta \bigr\}\subset {\cal G},
\eqno(\GopenA)$$
where the Morrey norm $\|\cdot\|_{M^{2,4/(p-1)}}$ is defined in (\DefMorrey).
In particular, there exists $\eta_1=\eta_1(u_0)>0$ such that
$$\bigl\{ v_0\in L^\infty(\Omega);\, 
\sup_{x\in\Omega}|x|^{2/(p-1)}|u_0(x)-v_0(x)| <\eta_1 \bigr\}\subset {\cal G}
\eqno(\GopenB)$$
and
$$\bigl\{ v_0\in L^\infty(\Omega);\, 
\|u_0-v_0\|_{n(p-1)/2} <\eta_1 \bigr\}\subset {\cal G}.
\eqno(\GopenC)$$
\smallskip
(ii) Let $\phi\in L^\infty(\Omega)$ and assume that $\lambda^*\in (0,\infty)$, where $\lambda^*$ is defined in (\DefLambdaStar).
If $\Omega$ is unbounded, assume in addition that $\phi\in BC^1(\Omega)$ and that $\phi$ satisfies (\hypdecayuzeto). 
Then $T_{max}(\lambda^*\phi)<\infty$.

In the above situation, it is a natural question whether or not $u^*$ can be continued 
after $T_{max}$ as some kind of weak solution.
In the case $\phi\ge 0$, then $\lambda\mapsto u(\cdot;\lambda \phi)$ is nondecreasing and
a suitable, unique minimal global weak continuation can be constructed 
by monotone approximation (cf.~[\NST]).\footnote{$^1$}{However, weak continuations need not be unique in general (see [\FM]).}
The behavior of the minimal global weak continuation for $p>p_S$ was studied in the radial case in [\GV]
and it was later shown [\CDZ] that if $\Omega$ is convex and bounded, then
$u^*$ becomes classical again for sufficiently large time and eventually decays uniformly to zero.
For sign-changing $\phi$, a global weak continuation was also constructed in [\CDZ] (see [\CDZ, Theorem~1(b))
by using compactness estimates and, as a consequence of the results in [\BS], the same behavior occurs.

Our next result extends this to the Cauchy problem (or to convex unbounded domains),
and also recovers the case of convex bounded domains with simpler proof. 
Actually, it applies to any pointwise limit of global classical solutions,
which in particular covers the weak continuations mentioned above.

\proclaim Theorem~5.
Assume $p>p_S$, $\Omega$ convex. Let $u$ be a pointwise limit of global classical solutions of (\equp), namely
$$\lim_{j\to\infty} u_j(x,t)=u(x,t),\quad\hbox{ a.e. in $\Omega\times(0,\infty)$,}$$
where the initial data of $u_j$ are of the form $u_{0,j}=\lambda_j\phi$ with $\phi\in L^\infty(\Omega)$
and $\lambda_j\to\lambda_0\in (0,\infty)$. If $\Omega$ is unbounded, assume in addition that $\phi\in BC^1(\Omega)$ and
that $\phi$ satisfies~(\hypdecayuzeto).
Then there exists $t_0>0$ such that $u$ is a bounded classical solution on $\Omega\times (t_0,\infty)$.
Moreover $u$ satisfies the decay property (\Sdecayu).

\medskip

The outline of the rest of the article is as follows.
Several remarks, especially about the optimality of the results, as well as ideas of proofs,
are given in Section~3. The preliminary Section~4 consists of two subsections of reminders, respectively on Morrey spaces
and on Giga-Kohn weighted energy.
In Sections~5 and 6, we then present our two essential technical results: Proposition~5.1
on the derivation of $M^{2,4/(p-1)}$-Morrey estimates from weighted energy bounds, 
and Proposition~6.1 on the global existence and decay for small initial data
 in the critical Morrey space $M^{2,4/(p-1)}(\Omega)$.
The main results are then deduced as consequences thereof in Section 7,
except for Theorem 4, which is proved in Section~8, as a consequence of a suitable continuous dependence property 
for problem (\equp) with respect to the critical Morrey norms (Proposition~8.1).
Finally, in appendix, for convenience and self-containedness, we provide the proofs of some important known results that we 
have used, concerning $L^p$-$L^q$ estimates for the linear heat semigroup
in Morrey spaces (Proposition~4.1),
basic properties of Giga-Kohn weighted energy (Proposition~4.3),
plus some auxiliary lemmas.

 \medskip\medskip
\goodbreak

{\bf 3. Discussion and ideas of proofs.}
\medskip

{\bf Remarks 3.1.}
(i) In the case $\Omega=\R^n$, unlike in the case of bounded domains, boundedness of global solutions
cannot  hold in general for supercritical $p$ under the mere assumption $u_0\in L^\infty$. 
Actually, the spatial decay assumption (\SdecayHba) in Theorem~1 is essentially optimal for the boundedness of global solutions,
Indeed, for all $p>p_S$, there exist {\it unbounded} global solutions for initial data 
with the borderline decay 
$$u_0(x)\sim |x|^{-2/(p-1)}\quad\hbox{as $|x|\to\infty$}$$ (cf.~[\PY, \PYb]).
On the other hand, the solution of (\equp) is never global if the spatial decay is slower, namely if $\displaystyle\lim_{|x|\to\infty}|x|^{2/(p-1)}u_0(x)=\infty$ (see~[\LN, \QSb]).
\smallskip

(ii) The conclusion (\Sdecayu) about temporal decay in Theorem~1, is also optimal.
Indeed (see~e.g. [\QSb, Chapter 20]), for any $p>(n+2)/n$
and $k>1/(p-1)$, if we take $u_0(x)=\eps(1+|x|^2)^{-k}$ with $\eps>0$ small (which satisfies (\SdecayHb)),
then $u$ is global and has the decay rate $\|u(t)\|_\infty\sim t^{-k}$, which can be made arbitrarily close to 
$t^{-1/(p-1)}$. Moreover, the conclusion (\Sdecayu) is in contrast with the subcritical case $p<p_S$ since in this case,
there exist initial data $u_0$ with exponential spatial decay, such that $u$ is global and 
$$\|u(t)\|_\infty\sim t^{-1/(p-1)} \quad\hbox{as $t\to\infty$},$$
including forward self-similar solutions (see [\HW, \Kaw]).
\smallskip

(iii) Under assumption (\SdecayHaaa) in Theorem 1 (or (\SdecayHaa) in Corollary~3), the solution $u$ actually enjoys
better decay properties, depending on the value of $q$.
In fact, $u$ behaves like the linear part of the equation,
in the sense that there exists a constant $C>0$ (depending on $u$)
such that 
$$|u(t)|\le Ce^{-tA}|u_0|\quad\hbox{ in $\Omega\times (0,\infty)$,}
\eqno(\LinearBehavior)$$
 where $e^{-tA}$ is the (Dirichlet) heat semigroup in $\Omega$.
In particular, under assumption (\SdecayHaaa) in Theorem 1, we have
$$\|u(t)\|_\infty\le Ct^{-(n-q)/2q},\quad t\ge 1,
\eqno(\LinearBehaviorB)$$
where $(n-q)/2q>1/(p-1)$;  see~Remark~7.1.

\smallskip

(iv) For $p>p_S$, $\Omega=\Rn$ and $0\le u_0\in L^\infty$ radially symmetric, the boundedness of global solutions
was already known in some cases. It is true if we assume $u_0\in  H^1(\R^n)$ (see~[\MMJFA, Theorem 5.15]).
It is also true if $p>p_{JL}$ and
$u_0(x)\leq (1-\eps)U_*(|x|)$ for $|x|$ large and some $\eps>0$,
where $p_{JL}= 1+4{n-4+2\sqrt{n-1}\over (n-2)(n-10)_+}$ and $U_*$ is the singular steady-state
(see~[\MizGlobal]),
or if $p<p_{JL}$ and $u$ satisfies certain finite intersection number properties 
with $U_*$ (see~[\Miz, Lemma 2.2]). 

\smallskip

(v) For $p=p_S$, the situation is quite different. Indeed if $\Omega$ is a ball, there exist
global positive, radial solutions such that $\lim_{t\to\infty}\|u(t)\|_\infty=\infty$ 
(see~[\GV]), and their grow-up rates are precisely known [\GK].
This is conjectured to happen also when $\Omega=\R^n$ if $n\le 4$ (see~[\FK]).
\fin

\smallskip\smallskip

{\bf Remarks~3.2.}
(i) When $\Omega$ is bounded, the rate (\Sdecayu)  can of course be improved. 
Indeed, in a bounded domain, any decaying solution of (\equp) decays at an exponential rate, 
since it is well known that the trivial solution is exponentially asymptotically stable.
This remains true if $\Omega$ is unbounded but has finite inradius; see~[\SouJFA].

\smallskip

(ii) In the subcritical case $p<p_S$, it is well known [\Gig, \Qui] that global solutions are not only bounded but satisfy an a priori
estimate of the form 
$$\sup_{t\ge 0}\|u(t)\|_\infty\le K(\|u_0\|_\infty),
\eqno(\AprioriEst)$$
 where the constant $K$ is bounded on bounded sets.
Such an estimate implies in particular that the borderline solutions (cf.~Section~2.2) are global and bounded.
Consequently, for $p>p_S$ in convex domains, estimate (\AprioriEst) cannot be true, in view of Theorem~4.
However, by the proof of Theorem~2 and Corollary~3, the following weaker, {\it delayed}
a priori estimate of global solutions is true:
$$\exists t_0=t_0(\|u_0\|)>0,\qquad \sup_{t\ge t_0}\|u(t)\|_\infty\le K(\|u_0\|).$$
Here $t_0$ is bounded on bounded sets and we may take $\|u_0\|=\|u_0\|_\infty$ if $\Omega$ is convex bounded, 
and e.g., $\|u_0\|=\|u_0\|_q+\|\nabla u_0\|_m$ with $q\in [2,{n(p-1)\over 2})$ and $m\in [2,{n(p-1)\over p+1})$
otherwise.
This estimate, which remains valid for borderline solutions, is consistent 
with the fact that any borderline solution blows up in a finite time $T^*$
(which thus satisfies $T^*<t_0(\|u_0\|)$).

 \smallskip

(iii) Theorem~1 and Corollary~3 remain true if, 
instead of (\SdecayHaaa) or (\SdecayHaa), we assume $u_0\in L^q(\Omega)$ for some $q\in [2,{n(p-1)\over p+1})$.
Indeed, it is easy to show that condition (\SdecayHaa) is then satisfied at positive time. 
This avoids any assumption on $\nabla u_0$
(which need not exist), but at the expense of a possibly more restrictive condition at infinity on $u_0$.
On the other hand, it should be possible to extend Theorem~1
(where the integrability or decay assumption is made on $\nabla u_0$ only) to general unbounded domains,
at the expense of significant additional technicalities. However, for the sake of simplicity, we have refrained ourselves from
treating this.

\smallskip

(iv) It is still an open problem whether or not boundedness of global solutions remains true without 
the convexity assumption on $\Omega$.
If $p\ge p_S$ and $\Omega$ is bounded and starshaped (in particular if it is convex), then any bounded global solution of
(\equp) has to decay in $L^\infty$, 
owing to the existence of a Liapunov functional and to the absence of nontrivial steady states
due to Pohozaev's identity.
But the decay of all global solutions cannot be true for general bounded domains, 
since positive stationary solutions exist for all $p>1$ when $\Omega$ is for instance an annulus.
\fin

\medskip\medskip

The original proofs of Theorem~A for bounded domains in [\BS] and [\CDZ] (and nonnegative solutions in the latter)
are based on two ingredients:
\goodbreak

\smallskip
\item{(a)} a monotonicity property, which yields suitably decaying space-time estimates on the solution;
\smallskip
\item{(b)} an $\varepsilon$-regularity result in terms of the space-time estimates.
\smallskip
Concerning (a), the authors of [\CDZ] use the nondecreasing property of the weighted energy of Giga and Kohn [\GKa, \GKb] for equation (\equp) rewritten in similarity variables, with rescaling time as a free parameter 
(see [\MMCPAM] for earlier related arguments in the radial framework).
This argument is combined with a Pohozaev-type identity in the original variables.
As for [\BS], a different tool is used, namely a parabolic analogue of an elliptic monotonicity formula obtained and used in [\Pac]
to prove partial regularity of weak solutions of the stationary equation $-\Delta u=|u|^{p-1}u$ for $p\ge p_S$.
Some of the close relations between these two points of view are noticed in [\BS]; see~[\BS, Section~3.1 and~p.~2348].

As for (b), the $\eps$-regularity result in [\CDZ]  is derived by a clever use of linear parabolic $L^q$-estimates,
leading to local uniform bounds of the solution. 
A different idea is used in [\BS], based on a novel analysis of
equation (\equp) in a parabolic Morrey space, involving some delicate bootstrap arguments. 

We shall here develop a modified approach, which also covers possibly unbounded domains
and gives the optimal decay rate in the case of $\Rn$.
It actually combines the Giga-Kohn weighted energy framework with an analysis of (\equp) in
an {\it elliptic} Morrey space of functions of $x$ (instead of a parabolic Morrey space of functions of $x, t$).
This new procedure has the advantage to replace the $\eps$-regularity issue with a 
smoothing and decay property for an initial value problem. 

Comparing with (a)(b) above, our strategy can be summarized 
in the following two steps:

\smallskip
\item{(a')} By using properties of the weighted energy of Giga and Kohn,
we show that for any global classical solution, 
the critical Morrey norm $\|u(t)\|_{M^{2,4/(p-1)}(\Omega)}$ decays 
as~$t\to\infty$;\footnote{$^1$}{This is where the condition $p>p_S$ enters 
and it also requires $u_0$ with suitable spatial decay in the case of unbounded domains}
 \smallskip
\item{(b')} We prove that any solution starting from an initial data with sufficiently small
$M^{2,4/(p-1)}$-norm exists globally and decays in $L^\infty$ as $t\to\infty$. 

\smallskip

Property (b') is proved by a semigroup argument, in the line of [\Wei, \SnTW, \SouCPDE]
which, as another important advantage, also yields the correct temporal decay of the solution in $L^\infty$
in the case $\Omega=\Rn$.
We here make use of linear smoothing estimates of Kato [\Kato] for the heat semigroup in Morrey spaces.

We note that our proof of Theorem A (bounded domains) is somewhat simpler than those of [\BS] and [\CDZ].
However the tools developed in [\BS, \CDZ] allow to obtain further delicate properties
such as estimates of the Hausdorff dimension of the singular set
of certain blowup solutions for $p>p_S$, which we do not consider.
\smallskip

\medskip

\goodbreak

{\bf 4. Preliminaries.}
\medskip
\ \ {\bf 4.1. Morrey spaces.}
\medskip

Let $\Omega$ be any domain in $\Rn$, $q\in [1,\infty)$ and $\lambda\in [0,n]$.
Recall that the Morrey space $M^{q,\lambda}(\Omega)$ is defined by 
$$M^{q,\lambda}(\Omega)=\bigl\{f \in L^q_{loc}(\Omega):\,
\|f\|_{M^{q,\lambda}}<\infty\bigr\}, 
$$
where
$$\|f\|^q_{M^{q,\lambda}}=\|f\|^q_{M^{q,\lambda}(\Omega)}:=\sup_{a\in\Omega}\  \sup_{R>0} r^{\lambda-n}\int_{B_R(a)\cap \Omega}|f|^q\, dx.
\eqno(\DefMorrey)$$
Also we set $M^{\infty,\lambda}(\Omega)=L^\infty(\Omega)$.
We note that one can use instead
$$\sup_{a\in\Rn}\  \sup_{R>0} r^{\lambda-n}\int_{B_R(a)\cap \Omega}|f|^q\, dx,$$
which is easily seen to give an equivalent norm.
Observe also that $M^{q,n}(\Omega)=L^q(\Omega)$, whereas $M^{q,0}(\Omega)=L^\infty(\Omega)$,
owing to the Lebesgue differentiation theorem.
For all $1<p\le r<\infty$, $m\in (0,r]$ and $\lambda\in [0,n]$, we have
$$\||f|^m\|_{M^{r/m,\lambda}}=\|f\|^m_{M^{r,\lambda}},
\qquad f\in M^{r,\lambda}(\Omega), \eqno(\PropMorreyA)$$
and, as a consequence of H\"older's inequality,
$$\|fg\|_{M^{r/p,\lambda}}\le\|f\|_{M^{r/(p-1),\lambda}} \, \|g\|_{M^{r,\lambda}},
\qquad f\in M^{r/(p-1),\lambda}(\Omega), \ g\in M^{r,\lambda}(\Omega). \eqno(\PropMorreyB)$$

\smallskip\smallskip

Let $e^{-tA}$ be the Dirichlet heat semigroup on $L^\infty(\Omega)$.
By results in [\Kato], it is known that $e^{-tA}$ enjoys good smoothing properties in the scale of Morrey spaces.

\proclaim Proposition~4.1. 
Let $\Omega$ be an arbitrary domain of $\Rn$, let $1\le p\le q\le \infty$ and $0\le\lambda\le n$.
Then there exists a constant $C=C(n,p,q,\lambda)>0$ such that,
for all $f\in L^\infty(\Omega)\cap M^{p,\lambda}(\Omega)$, 
$$\|e^{-tA}f\|_{M^{q,\lambda}}\le Ct^{-{\lambda\over 2}({1\over p}-{1\over q})}\|f\|_{M^{p,\lambda}},
\quad t>0.
\eqno(\RegulMorrey)$$
Moreover, 
$$\|e^{-tA}f\|_{M^{p,\lambda}}\le \|f\|_{M^{p,\lambda}},\quad t>0.
\eqno(\ContractionMorrey)$$

\smallskip
A proof is provided in Appendix for convenience.
\medskip

{\bf Remark 4.1.} 
Observe that in the case $\lambda=n$, since $M^{p,n}(\Omega)=L^p(\Omega)$,
this corresponds to the usual $L^p$-$L^q$ estimate,
and that the case $\lambda=0$ is just $\|e^{-tA}f\|_\infty\le C \|f\|_\infty$.
\fin
\medskip

We next note the following elementary properties which will be quite useful,
since they relate various ``critical'' norms associated with problem (\equp).
Here $L^\infty_{2/(p-1)}(\Omega)$ denotes the weighted $L^\infty$ space with weight $|x|^{2/(p-1)}$,
with norm
$$\|f\|_{\infty, 2/(p-1)}:={\rm ess}\sup_{x\in\Omega} |x|^{2/(p-1)}|f(x)|.
\eqno(\WeightedLinfty)$$

\goodbreak

\proclaim Proposition~4.2.
Let $\Omega$ be an arbitrary domain of $\Rn$ and $p>1+{2\over n}$.
\smallskip
(i) For all $1\le q<n(p-1)/2$, we have
$$L^{n(p-1)/2}(\Omega) \hookrightarrow M^{q,2q/(p-1)}(\Omega)
\eqno(\compMorreyGnormA)$$
and
$$L^\infty_{2/(p-1)}(\Omega)\hookrightarrow M^{q,2q/(p-1)}(\Omega).
\eqno(\compMorreyGnormB)$$

\smallskip
(ii) For all $q\in [1,\infty)$, $\lambda\in [0,n]$ and $\phi\in L^\infty(\Omega)$, we have
 $$\|\phi\|_{M^{q,\lambda}}\le C \sup_{t>0}t^{\lambda/2}\|G_t*|\phi|^q\|_\infty.
\eqno(\compMorreyGnorm)$$
In particular, for $p>1+{4\over n}$, if $u_0\in BC^1(\Omega)$ satisfies (\hypdecayuzeto), then
$u_0\in M^{2,4/(p-1)}(\Omega)$. 

\medskip

{\bf Remark 4.2.} 
Although we shall not use this property, let us mention that, for $\lambda>0$, the quantity 
$\sup_{t>0}t^{\lambda/2}\|G_t*f\|_\infty$
is known to be equivalent to the norm of the homogeneous Besov space $\dot B^{-\lambda}_{\infty,\infty}$
(cf.~[\Tri, p.~192]) and that condition (\hypdecayuzeto) can thus be interpreted in terms of 
belonging to suitable closed subspaces of Besov spaces.
\fin
 \medskip

{\bf Proof.} Let $m=n(p-1)/2$. For all $a\in\Omega$, $R>0$, by H\"older's inequality, we have
$$R^{{2q\over p-1}-n}\int_{B_R(a)\cap \Omega}|f|^q\,dx
\le CR^{{2q\over p-1}-n}R^{n(1-{q\over m})}\Bigl(\int_{B_R(a)\cap \Omega}|f|^m\,dx\Bigr)^{q\over m}
\le C\|f\|_m^q$$
and (\compMorreyGnormA) follows.
\smallskip

To prove (\compMorreyGnormB), let us show that 
$$|x|^{-2/(p-1)}\in M^{q,\lambda}(\Omega),\quad \lambda=2q/(p-1).
\eqno(\compMorreyintermed)$$
Setting $k=2/(p-1)$, we write
$$R^{{2q\over p-1}-n}\int_{B_R(a)\cap\Omega}|x|^{-kq}\,dx
\le R^{kq} \int_{B_1} |a+Ry|^{-kq}\,dy
= \int_{B_1} |aR^{-1}+y|^{-kq}\,dy.$$
If $aR^{-1}\ge 2$, then 
$\int_{B_1} |aR^{-1}+y|^{-kq}\,dy\le C:=\int_{B_1} \,dy.$
If $aR^{-1}<2$, since $kq<n$, we have
$$\int_{B_1} |aR^{-1}+y|^{-kq}\,dy=\int_{B_1(aR^{-1})} |y|^{-kq}\,dy\le C:=\int_{B_3} |y|^{-kq}\,dy.$$
This proves (\compMorreyintermed). For all $f\in L^\infty_{2/(p-1)}(\Omega)$, since
$|f|\le \|f\|_{\infty, 2/(p-1)} |x|^{-2/(p-1)}$, it follows that $f\in M^{q,\lambda}(\Omega)$ 
and $\|f\|_{M^{q,\lambda}} \le \||x|^{-2/(p-1)}\|_{M^{q,\lambda}}\|f\|_{\infty, 2/(p-1)}$,
and (\compMorreyGnormB) is proved.
\smallskip

Finally, assertion (\compMorreyGnorm) is a consequence of
$$R^{\lambda-n}\int_{B_R(a)\cap\Omega}|\phi|^q\,dx 
\le CR^{\lambda-n}\int_{B_R(a)\cap\Omega} R^n G_{R^2}(x-a)|\phi(x)|^q\,dx
\le CR^\lambda\|G_{R^2}*|\phi|^q\|_\infty.$$
\fin

\medskip
\goodbreak

\ \ {\bf 4.2. Similarity variables and weighted energy.}
\medskip

We here recall some fundamental tools from Giga and Kohn [\GKa, \GKb]. 
Let $\Omega$ be a (possibly unbounded) domain of $\Rn$ of class $C^{2+\alpha}$.
Given $a\in\overline\Omega$, 
the backward similarity variables $(y,s)$ with respect to $(a,T)$ are given by
$$ y:={x-a\over \sqrt{T-t}},\qquad s:=-\log(T-t).
\eqno(\SGKdefSelfsimilA)
$$
Let $p>1$ and set $\beta={1\over p-1}$. 
If $u$ is a solution of (\equp) on $(0,T)$, then the
corresponding rescaled solution $w=w_{a,T}(y,s)$ is given by
$$ w(y,s)=w_{a,T}(y,s):=e^{-\beta s} u(a+e^{-s/2}y,T-e^{-s}),
\eqno(\SGKdefSelfsimilB)
$$
and it is defined for all $s>s_0:=-\log T$ in the rescaled domain 
$$D(s):=\bigl\{y\in\R^n;\ a+e^{-s/2}y\in \Omega\bigr\}=e^{s/2}(\Omega-a).
\eqno(\SGKdefrescaledOmega)
$$
By direct calculation, $w$ satisfies the equation
$$ \rho w_s-\nabla\cdot(\rho\nabla w)=\rho |w|^{p-1}w-\beta\rho w,\qquad \rho(y)=e^{-{|y|^2\over 4}},
\eqno(\QGKeqrho)
$$
in $\bigcup_{s>s_0}D(s)\times\{s\}$. 
Finally, the weighted energy $E(s)=E(w_{a,T}(s))$ is given~by
$$ E(w_{a,T}(s)):=\int_{D(s)}\Bigl({1\over 2}|\nabla w|^2+{\beta\over 2} w^2
       -{1\over p+1}|w|^{p+1}\Bigr)\rho\,dy.
\eqno(\SGKdefenergyOmega)
$$
The main properties of the rescaled energy $E$, in particular its time monotonicity,
 are described in the following proposition [\GKa, \GKb].

\proclaim Proposition~4.3.
Let $a\in\Omega$ and assume that $\Omega$ is starshaped with respect to~$a$.
Let $p>1$, $u_0\in BC^1(\Omega)$ and assume that the solution $u$ of (\equp) exists on $(0,T)$.
Then, under the notation of the above paragraph, for all $s>s_0$, we have
$$
{1\over 2}{d\over ds}\int_{\Rn}w^2\rho\,dy
=-2E(w(s))+{p-1\over p+1}\int_{\Rn}|w|^{p+1}\rho\,dy,
\eqno(\QGKidA)
$$
$$
{d\over ds}E(w(s)) \le -\int_{\Rn}w_s^2\rho\,dy,
\eqno(\QGKidB) 
$$
$$E(w(s))\geq 0
\eqno(\QGKenergypos) $$
and
$$\int_{\Rn}w^2\rho\,dy\leq C(n,p) \bigl[E(w(s_0))\bigr]^{2/(p+1)}.
\eqno(\QGKLtwocontrol)
$$

The proof is provided in Appendix for convenience.

\medskip\medskip
{\bf 5. Estimate of $u(t)$ in the Morrey space $M^{2,4/(p-1)}(\Omega)$.}
\medskip

Our first main technical result is the following proposition, which enables one
to relate the Morrey norm $\|u(t)\|_{M^{2,4/(p-1)}(\Omega)}$ to weighted energies 
for appropriate rescaling times~$T>t$. The latter can be estimated 
by suitable norms of the initial data, owing to the monotonicity of the weighted energy. 
This will eventually lead to the decay of $\|u(t)\|_{M^{2,4/(p-1)}(\Omega)}$. 

\proclaim Proposition~5.1. 
Let $\Omega\subset\Rn$ be a (possibly unbounded) domain of class $C^{2+\alpha}$, 
let $u_0\in BC^1(\Omega)$ and let $u$ be the solution of~(\equp).
\smallskip
(i) Let $T>0$, $a\in\overline\Omega$ and assume that $\Omega$ is starshaped with respect to $a$.
If $u$ exists on $(0,T)$, then for all $t\in (0,T)$, we have
$$(T-t)^{{2\over p-1}-{n\over 2}}\int_\Omega e^{-{|x-a|^2\over 4(T-t)}} u^2(x,t)\,dx
\le C(n,p){\cal A}^{2\over p+1}(T,u_0,a),$$
where
$${\cal A}(T,u_0,a)=\Bigl[T^{p+1\over p-1}G_T*(|\nabla u_0|^2)+T^{2\over p-1}G_T*(|u_0|^2)\Bigr](a).$$
 \smallskip
(ii) Assume $\Omega$ convex and set $\mu={4\over p-1}$. 
If $u$ is global, then, for all $t_0>0$, we have
$$\|u(t_0)\|_{M^{2,\mu}(\Omega)} \le C(n,p){\cal N}^{1\over p+1}(u_0,t_0),$$
where
$${\cal N}(u_0,t_0)
=\sup_{t\ge t_0} \Bigl(t^{p+1\over p-1}\|G_t*(|\nabla u_0|^2)\|_\infty+t^{2\over p-1}\|G_t*(|u_0|^2)\|_\infty\Bigr).$$

{\bf Proof.}
In this proof, $C$ denotes a generic positive constant depending only on $n$ and $p$.
\smallskip

(i) Let $(y,s)$ be the backward similarity variables with respect to $(a,T)$, 
$w=w_{a,T}(y,s)$ the corresponding rescaled solution,
$D(s)$ the rescaled domain, and $E(s)=E(w_{a,T}(s))$ the weighted energy 
(cf.~(\SGKdefSelfsimilA)--(\SGKdefenergyOmega)). 
By Proposition~4.3, for all $s\ge s_0=-\log T$, we have
$$0\le E(s)\le E(s_0),\qquad \int_{D(s)} |w|^2 \rho\, dy\le C[E(s_0)]^{2/(p+1)}.
$$
Switching back to the original variables, it follows that
$$(T-t)^{{2\over p-1}-{n\over 2}} \int_\Omega e^{-{|x-a|^2\over 4(T-t)}}u^2(x,t)\, dx
\le C[E(s_0)]^{2/(p+1)}.
$$
To estimate the RHS, we write
$$  \int_{D(s_0)}|\nabla w(y,s_0)|^2\rho\, dy
=T^{{p+1\over p-1}-{n\over 2}} \int_\Omega |\nabla u_0|^2\,e^{-{|x-a|^2\over 4T}}\, dx
=CT^{p+1\over p-1}(G_T*|\nabla u_0|^2)(a)
\eqno(\RHSa)$$
and
$$
\int_{D(s_0)} w^2(y,s_0)\rho\, dy
=T^{{2\over p-1}-{n\over 2}} \int_\Omega |u_0|^2\,e^{-{|x-a|^2\over 4T}}\, dx
=CT^{2\over p-1}(G_T*|u_0|^2)(a),
\eqno(\RHSb)$$
Since 
$$E(w_{a,T}(s))\le \int_{D(s)}\Bigl({1\over 2}|\nabla w|^2+{\beta\over 2} w^2\Bigr)\rho\,dy,
$$
this guarantees the assertion.

\medskip

(ii) Pick $t_0>0$ and $a\in\Omega$. Since $\Omega$ is convex, it is starshaped with respect to $a$.
For any $R>0$, we choose $T=t_0+R^2$. For each $x\in B(a,R)$, we have ${|x-a|^2\over T-t_0}\le 1$.
It then follows from assertion (i) that
$$\eqalign{
R^{\mu-n} \int_{\Omega\cap B(a,R)} u^2(x,t_0)\,dx
&\le C(T-t_0)^{\mu-n\over 2}\int_\Omega e^{-{|x-a|^2\over 4(T-t_0)}} u^2(x,t_0)\,dx \cr
&\le  C{\cal A}^{2\over p+1}(T,u_0,a)\le C{\cal N}^{2\over p+1}(t_0,u_0).
}$$
The assertion follows by taking supremum over $a\in\Omega$ and $R>0$. \fin

\medskip
{\bf Remark~5.1.} Assume $u_0\in H^1(\Omega)$ (which is always true in bounded domains, 
starting from $L^\infty$, after a time shift).
Then, since $p>p_S$ is equivalent to ${p+1\over p-1}-{n\over 2}<0$, it is already clear from formulae (\RHSa)-(\RHSb) that
$$E(s_0)=E\bigl(w_{a,T}(-\log T)\bigr)\to 0,\quad T\to\infty.$$
This crucial observation is a starting point of the analysis of [\BS]
(although, unlike in our case, it is used there in conjunction with space-time estimates and parabolic Morrey norms).
\fin

\medskip
\medskip
{\bf 6. Global existence and decay for small data in critical Morrey spaces.}
\medskip

The remaining task is now to infer the uniform decay of $u(t)$ from its decay in the Morrey space
$M^{2,4/(p-1)}$. This space --~and more generally $M^{q,2q/(p-1)}$ for $p\ge 1+{2q\over n}$~--
turns out (when $\Omega=\Rn$) to be invariant by the scaling of the equation, namely 
$$\lambda\mapsto u_\lambda(x,t)=\lambda^{2/(p-1)}u(\lambda^2t,\lambda x).
\eqno(\DefScaling)$$
This property is similar to that of the critical $L^q$ space with $q=q_c=n(p-1)/2$ 
(see, e.g., [\QSb, Chapter 20] and also [\Kar] for further examples).
Owing to the linear smoothing estimates of [\Kato], the following small data global existence and decay result for problem (\equp) 
 can be shown in a similar way as for $L^{q_c}$,
based on ideas from [\Wei, \SouCPDE, \SnTW].
In the case of Morrey spaces, results of this type were first obtained in [\GMi, \Kato] for the Navier-Stokes system.
For recent developments concerning other problems such as chemotaxis systems, see e.g. [\LR, \BKZ].
Related results for semilinear parabolic equations appear in [\YZ, \Tay, \Kar], but they do not seem suitable to our needs.

\proclaim Proposition~6.1. 
Let $\Omega\subset\Rn$ be a (possibly unbounded) domain of class $C^{2+\alpha}$. 
Let $q\in (1,\infty)$ and $p\ge 1+{2q\over n}$.
There exist $\eps_0=\eps_0(n,p,q)$ and $C_0=C_0(n,p,q)>0$ with the following property.
For any $u_0\in L^\infty(\Omega)\cap M^{q,2q/(p-1)}(\Omega)$, if 
$$\|u_0\|_{M^{q,2q/(p-1)}(\Omega)}\le \eps_0,$$
then the corresponding solution $u$ of (\equp) is global and satisfies
$$\sup_{t>0} t^{1\over p-1}\|u(t)\|_\infty\le C_0 \|u_0\|_{M^{q,2q/(p-1)}(\Omega)}.
$$

\goodbreak
\medskip
{\bf Remarks~6.1.} 
(i) Although this is out of the scope of this paper, we mention that local well-posedness for problem (\equp)
is not expected to hold if we only assume $u_0\in M^{q,2q/(p-1)}(\Omega)$ (cf.~[\Kato, \Tay] for related issues 
concerning the Navier-Stokes system; this is related with the lack of density of smooth functions in Morrey spaces). 
Actually, locall well-posedness should be true under the slightly stronger assumption
$$u_0\in \hat M^{q,2q/(p-1)}(\Omega)=\Bigl\{u_0\in M^{q,2q/(p-1)}(\Omega);\, 
\lim_{r\to 0} \ \sup_{a\in\Omega}\  r^{\lambda-n}\int_{B_R(a)\cap \Omega}|f|^q\, dx =0
\Bigr\}.$$

(ii) Proposition~5.1 guarantees that any borderline global weak solution $u^*$
satisfies 
$$\sup_{t>0}\|u^*(t)\|_{M^{2,4/(p-1)}(\Omega)}<\infty.
\eqno(\SupMorrey)$$
 In spite of this, $u^*$ blows up in $L^\infty$ norm 
at some finite time $T_{max}(\lambda^*\phi)$ by Theorem~4.
This implies in particular that the (classical) existence time of a solution 
is not uniform for initial data in bounded sets of $M^{2,4/(p-1)}$.
Likewise, if $u_0$ has the borderline decay 
(i.e., $O$~instead of $o$ in (\SdecayHb)), then any global classical solution satisfies (\SupMorrey),
by Proposition~5.1 and the proof of Corollary~5.
However global unbounded solutions do exist for such $u_0$ (cf.~Remark~3.1(i)).

\smallskip
(iii) The above nonuniformity phenomenon is typical when dealing with critical spaces
and it is known to occur for instance for problem (\equp) in the critical Lebesgue space $L^q(\Omega)$ 
with $q=n(p-1)/2$ (cf.~[\QSb, Sections 15-16]).
However, $L^q(\Omega)$ is a strict subspace of $M^{2,4/(p-1)}(\Omega)$
with stronger norm (cf.~Proposition~4.2)
and, as an interesting difference, the critical Morrey norm here does not blow up at $T_{max}$,
whereas, in all the known examples, the critical $L^q$ norm always does
(although the question is still open in general).
Since, on the other hand, problem (\equp) is locally well posed in the critical $L^q$ space,
this difference may be also connected with remark~(i).
Related blowup results in the critical Lorentz space $L^{q,\infty}$ can be found in [\GKb]  for $p<p_S$.

\smallskip
(iv) Consider the space $L^\infty_{2/(p-1)}(\Omega)$, defined in (\WeightedLinfty),
which is also critical (in the sense that it is invariant by the scaling 
(\DefScaling)),
and is another strict subspace of $M^{2,4/(p-1)}(\Omega)$ (cf.~Proposition~4.2).
It is known [\MMJFA] that in the radial case, for $p>p_{JL}$, the $L^\infty_{2/(p-1)}$ norm may or may not blow up at $T_{max}$,
depending on the type of blowup (type I or II, unfocused or focused at $x=0$).
\fin
\bigskip

{\bf Proof of Proposition~6.1.}
Denoting by $e^{-tA}$ the Dirichlet heat semigroup on $L^\infty(\Omega)$ and setting $u^p=|u|^{p-1}u$,
we have, for all $t\in (0,T)$,
$$u(t)=e^{-tA}u_0+\int_0^t e^{-(t-s)A}u^p(s)\, ds$$
(in the $L^\infty$ sense),
by the variation-of-constants formula. Set $\lambda=2q/(p-1)\in (0,n]$ and denote $|\cdot |_m=\|\cdot\|_{M^{m,\lambda}(\Omega)}$ for all $m\in [1,\infty]$. 
Fix $r\in\bigl(\max(p,q),pq\bigr)$, so that
$$1<r/p<q<r$$
and set 
$$\beta={\lambda\over 2}\Bigl({1\over q}-{1\over r}\Bigr)={1\over p-1}\Bigl(1-{q\over r}\Bigr)<{1\over p}.$$
In all the proof, $C_i$ (resp., $C$) denote fixed (resp., generic) positive constants depending only on $n,p,q,r$.
By Proposition~4.1, we have 
$$|e^{-tA} \phi|_r\le C_1t^{-\beta}|\phi|_q,\qquad \phi\in M^{q,\lambda},\quad t>0
\eqno(\SMorreyInvEstA)$$
and
$$|e^{-tA} \phi|_r\le C_1t^{-{\lambda(p-1)\over 2r}}|\phi|_{r/p},\qquad \phi\in M^{r/p,\lambda},\quad t>0.
\eqno(\SMorreyInvEstAA)$$

For any given $t_0\in (0,T_{max})$, set $M=M(t_0)=\sup_{t\in [0,t_0]}\|u(t)\|_\infty^{p-1}<\infty$.
On the interval $(0,t_0)$, we get $|u_t-\Delta u|\le M|u|$ hence, by the maximum principle, 
$$|u(t)|\le e^{Mt}e^{-tA}|u_0|\quad\hbox{ in $\Omega\times [0,t_0]$.}
\eqno(\SMorreyInvEstAB)$$
In particular $t^\beta|u(t)|_r\le 2C_1|u_0|_q$ for $t>0$ small.
We may thus define 
$$\tau_1=\sup\,\bigl\{t\in (0,T_{max});\,s^\beta|u(s)|_r\le 2C_1|u_0|_q\ \hbox{for all $s\in (0,t)$}\bigr\}$$
and we have $\tau_1\in (0,T_{max}]$. 
For all $0<s<t<\tau_1$, by (\SMorreyInvEstAA), (\PropMorreyA), we obtain 
$$
|e^{-(t-s)A}u^p(s)|_r
\leq C_1(t-s)^{-{\lambda(p-1)\over 2r}} |u^p(s)|_{r\over p}
=C_1(t-s)^{-{q\over r}} |u(s)|_r^p\le C|u_0|_q^p(t-s)^{-{q\over r}} s^{-\beta p}.
$$
On the other hand, since $1-\beta(p-1)-q/r=0$, we have
$$ 
t^\beta\int_0^t (t-s)^{-q/r} s^{-\beta p}\, ds
=t^{1-\beta(p-1)-q/r}\int_0^1 (1-\sigma)^{-q/r} \sigma^{-\beta p} \, d\sigma=C,\quad t>0,
\eqno(\intinv)$$
where the integrals are finite, due to $q/r<1$ and $\beta p<1$. 
It then follows that, for all $t\in (0,\tau_1)$,
$$t^\beta |u(t)|_r\leq t^\beta |e^{-tA} u_0|_r
+t^\beta \int_0^t |e^{-(t-s)A}u^p(s)|_r\, ds
\leq C_1|u_0|_q+C|u_0|_q^p.
\eqno(\intinvB)$$
Assume for contradiction that $\tau_1<T_{max}$. By continuity we may take $t=\tau_1$ in (\intinvB) to get 
$2C_1|u_0|_q\le C_1|u_0|_q+C|u_0|_q^p$, which is a contradiction if $|u_0|_q\le\eps_0$
with $\eps_0=\eps_0(n,p,q,r)>0$ sufficiently small. It follows that $\tau_1=T_{max}$ i.e.,
$$|u(t)|_r\le 2C_1|u_0|_q t^{-\beta},\quad 0<t<T_{max}.
\eqno(\SMorreyInvEst)
$$

Next, by Proposition~4.1, we have 
$$\|e^{-tA} \phi\|_\infty\le C_2t^{-\lambda/2r}|\phi|_r, \qquad \phi\in M^{r,\lambda},\quad t>0.
\eqno(\SMorreyInvEstB)$$
Let $C_0>0$ to be fixed later and set 
$$\tau_2=\sup\,\bigl\{t\in (0,T_{max});\,s^{1/(p-1)} \|u(s)\|_\infty \le C_0|u_0|_q\ \hbox{for all $s\in (0,t)$}\bigr\}\in (0,T_{max}].$$
Note that
$$u(t)=e^{-(t/2)A} u(t/2)
+\int_{t/2}^t e^{-(t-s)A} u^p(s)\, ds.
\eqno(\Varconsthalf)$$
Combining (\SMorreyInvEstB), (\SMorreyInvEst) and recalling $\beta={1\over p-1}-{\lambda\over 2r}$, it follows that
$$\eqalign{  
t^{1\over p-1}\|u(t)\|_\infty 
&\leq t^{1\over p-1} \|e^{-(t/2)A} u(t/2)\|_\infty
+t^{1\over p-1}\displaystyle\int_{t/2}^t \bigl\|u(s)\|^p_\infty\, ds\cr
&\leq 2^{\lambda\over 2r}C_2 t^{{1\over p-1}-{\lambda\over 2r}} |u(t/2)|_r+2^{1\over p-1}C_0^p |u_0|_q^p\cr
\noalign{\vskip 1mm}
&\leq 2^{{\lambda\over 2r}+\beta+1}C_1C_2 |u_0|_q+2^{1\over p-1}C_0^p |u_0|_q^p.
}$$
Choose $C_0=2^{{\lambda\over 2r}+\beta+2}C_1C_2$. 
Arguing as for $\tau_1$, taking $\eps_0$ smaller if necessary, we obtain $\tau_2=T_{max}$,
hence $T_{max}=\infty$
and the proposition follows. \fin

\medskip
{\bf 7. Proofs of Theorems 1, 2, 5 and Corollary 3.}
\medskip

{\bf Proof of Theorem~2.}
We may always assume that $u_0\in BC^1(\Omega)$. 
Indeed, if $\Omega$ is boun\-ded and $u_0\in L^\infty$ then
this is true after a time shift. Moreover, assumption (\hypdecayuzeto) is automatically satisfied if $\Omega$ is boun\-ded.
Indeed, since $|u_0|^2$ and $|\nabla u_0|^2$ then belong to $L^1(\Rn)$
(recall that they are extended by $0$ outside $\Omega$), we have
$$\|G_t*|\nabla u_0|^2\|_\infty+\|G_t*|u_0|^2\|_\infty\le Ct^{-n/2},$$
so that (\hypdecayuzeto) is true since ${p+1\over p-1}<{n\over 2}$ 
(due to $p>p_S$).\footnote{$^1$}{We observe in turn that, since $\|G_t*\phi\|_\infty\ge ct^{-n/2}$ for any nontrivial $\phi\ge 0$,
assumption (\hypdecayuzeto)  can be realized by a nontrivial $u_0$ {\it only} if $p>p_S$
(for any domain $\Omega$).}

Now, by assumption (\hypdecayuzeto) and Proposition~5.1, we have  
$$\displaystyle\lim_{t_0\to \infty}\|u(t_0)\|_{M^{2,4/(p-1)}(\Omega)}=0.$$
Pick any $\eps\in (0,\eps_0)$, where $\eps_0$ is given by Proposition~6.1 with $q=2$.
Then there exists $t_0=t_0(\eps)>0$ such that $\|u(t_0)\|_{M^{2,4/(p-1)}(\Omega)}\le \eps$,
and we deduce from Proposition~6.1 that
$$\|u(t)\|_\infty\le C_0\eps (t-t_0)^{-1/(p-1)}\le 2^{1/(p-1)}C_0\,\eps\, t^{-1/(p-1)},
\quad t\ge 2t_0.
\eqno(\Boundtzero)$$
Theorem 2 follows.  \fin

\goodbreak
\medskip

{\bf Proof of Theorem~5.}
By our assumptions, each $u_j$ satisfies
estimate (\Boundtzero) with uniform constants and a common $t_0$.
Passing to the limit, we see that $u$ satisfies the same estimate for a.e.~$t\ge 2t_0$.
By parabolic regularity, $u$ is in particular a classical solution
for $t\ge 2t_0$.
\fin

\smallskip

{\bf Proof of Corollary~3.}
Let us first assume (\SdecayHaa). We have $|u_0|^2\in L^{r}(\Omega)$ with $1\le r<n(p-1)/4$, hence
$$t^{2\over p-1}\|G_t*|u_0|^2\|_\infty\le C\| |u_0|^2\|_{r} \,t^{{2\over p-1}-{n\over 2r}}\to 0,\ t\to\infty.$$ 
Likewise, $|\nabla u_0|^2\in L^{m}(\Omega)$ with $1\le m<n(p-1)/2(p+1)$, hence
$$t^{p+1\over p-1}\|G_t*|\nabla u_0|^2\|_\infty\le C\| |\nabla u_0|^2\|_q Ct^{{p+1\over p-1}-{n\over 2m}}\to 0,\ t\to\infty.$$ 

Let us next assume (\SdecayHb). Set $k=2/(p-1)$ and choose $\delta>0$ such that $k+\delta<n/2$. 
By (\SdecayHb), for any $\eta>0$, there exists $A>0$ such that
$$|u_0(x)|^2\le \eta(1+|x|^2)^{-k}+A(1+|x|^2)^{-k-\delta}.$$
It is well-known that $\|G_t\ast (1+|x|^2)^{-k}\|_\infty \le C(t+1)^{-k}$ whenever $k\in (0,n/2)$
(see~e.g. [\QSb, Lemma 20.8]). Consequently,
$$t^{2\over p-1}\|G_t*|u_0|^2\|_\infty\le Ct^{2\over p-1}[\eta(t+1)^{-k}+A(t+1)^{-k-\delta}]\le 2C\eta$$
for $t$ large enough, hence $t^{2\over p-1}\|G_t*|u_0|^2\|_\infty\to 0$ as $t\to\infty$.
Similarly, we show that $t^{p+1\over p-1}\|G_t*|\nabla u_0|^2\|_\infty\to 0$ as $t\to\infty$.
\smallskip

In both cases, the conclusion then follows from Theorem~2.
\fin

We shall now derive Theorem~1 as a consequence of Corollary~3.
However, unlike in Corollary~3, 
$u_0$ is not assumed to satisfy any integrability or decay hypothesis (only $\nabla u_0$ does),
nor to belong to $BC^1$ in case of assumption (\SdecayHaaa).

To circumvent this, we shall use the following two general lemmas.
The first one allows to deduce a suitable integrability or decay property on a function 
from such an assumption on its gradient. The second one allows to perform a time shift,
so as to gain $BC^1$ regularity, while preserving assumption (\SdecayHaaa).

\bigskip

\proclaim Lemma~7.1. Let $n\ge 2$.
\smallskip
(i) Let $q\in [1,n)$ and let $v\in L^\infty(\Rn)$ be such that $\nabla v\in L^q(\Rn)$.
Then there exists a constant $K\in\R$ such that $v-K\in L^{q^*}(\Rn)$,
where $q^*=nq/(n-q)$.
\smallskip
(ii) Let $\beta>0$ and let $v\in BC^1(\Rn)$ be such that $|\nabla v|=o(|x|^{-\beta-1})$ as $|x|\to\infty$.
Then there exists a constant $K\in\R$ such that
$v(x)-K=o(|x|^{-\beta})$ as $|x|\to\infty$.

\medskip

\proclaim Lemma~7.2. Let $\Omega=\Rn$, $u_0\in L^\infty(\Rn)$ and let $u$ be the maximal 
classical solution of~(\equp). There exists $\tau\in (0,T_{max})$ such that
$|\nabla u(t)|\le 2G_t*|\nabla u_0|$ 
for all $t\in (0,\tau]$.

\medskip
\goodbreak

Results related with Lemma~7.1(i) can be found in, e.g., [\KS, Proposition~3.4] and [\OS].
Since they do not quite fulfill our needs, we give an elementary proof in Appendix.
We stress that in Lemma~7.1(i) we are not assuming $v$ to belong to the Sobolev space $W^{1,q}(\Rn)$
(otherwise, the conclusion is immediate by Sobolev embedding and necessarily $K=0$).
The homogeneous space $\dot W^{1,q}(\Rn)$, defined by completion of $C^\infty_0(\Rn)$ by the semi-norm $\|\nabla v\|_q$,
is sometimes used in similar situations. Note that addition of a constant cannot be a priori ruled out in this way.

As for Lemma~7.2, which is rather standard,
we also postpone its proof to the appendix, in order not to interrupt the main line of argument.

\medskip

{\bf Proof of Theorem~1.}
We first observe that we can always assume $u_0\in BC^1(\Rn)$,
making a time shift if necessary. Indeed, in case of assumption (\SdecayHaaa),
this assumption is still satisfied after replacing $u_0$ with $u(\tau)\in BC^1(\Rn)$ (for some small $\tau>0$),
thanks to Lemma~7.2.

Under assumption (\SdecayHaaa), by Lemma~7.1(i),
there exists $\psi\in L^{q^*}$ and  a constant $K\in\R$ such that $u_0=K+\psi$.
Since $q<n(p-1)/(p+1)$ implies $q^*<n(p-1)/2$, and $u_0\in BC^1$, it follows that
$$|u_0-K|^{p+1}+|\nabla u_0|^2\in L^m(\Rn)\quad\hbox{for some 
$m\in\bigl[1,{n(p-1)\over 2(p+1)}\bigr)$.}
\eqno(\TEnegA)$$
Under assumption (\SdecayHba), by Lemma~7.1(ii),
we may write $u_0=K+\psi$ with 
$\psi=o(|x|^{-2/(p-1)})$ as $|x|\to\infty$, hence
$$|u_0(x)-K|+|x|\,|\nabla u_0(x)|=o\bigl(|x|^{-{2\over p-1}}\bigr)\quad\hbox{as $|x|\to\infty$.}
\eqno(\TEnegB)$$

We claim that $T_{max}=\infty$ necessarily imposes $K=0$.
To prove this, we use the rescaled solution by similarity variables at $a=0$, $t=T$ and the corresponding weighted energy $E(s)=E(w_{0,T}(s))$, cf.~Section~4.2 
(a more direct comparison argument would apply in the case $u_0\ge 0$).
Assume for contradiction that $K>0$ (without loss of generality).
Denoting $\int=\int_{\Rn}$ and recalling the notation $G_T(x)=(4\pi T)^{-n/2}e^{-{|x|^2\over 4T}}$,
we have
$$\eqalignno{
T^{-{p+1\over p-1}}
&E(w_{0,T}(-\log T)) \cr
&=C_1\int |\nabla\psi|^2\,G_T+{C_2\over T} \int |K+\psi|^2\,G_T
-C_3\int |K+\psi|^{p+1}\,G_T \cr
&\le C_1\int |\nabla\psi|^2\,G_T
+{C_2\over T} \Bigl(1+ \int |K+\psi|^{p+1}\,G_T\Bigr)
-C_3\int |K+\psi|^{p+1}\,G_T\cr
&\le C_1\int |\nabla\psi|^2\,G_T
+{C_2\over T}-{C_3\over 2}\int |K+\psi|^{p+1}\,G_T\cr
&\le C\int (|\nabla\psi|^2+|\psi|^{p+1})\,G_T
+{C_2\over T}-C_4 K^{p+1}, &(\TEneg)}$$
for $T>0$ large enough, where we used $\int G_T=1$,
and also $|K+\psi|^{p+1}\ge 2^{-p} K^{p+1}-|\psi|^{p+1}$ in the last inequality.
On the other hand, by either (\TEnegA) or (\TEnegB), 
we have $|\nabla\psi|^2+|\psi|^{p+1}\in L^k(\Rn)$ for some $k\in(1,\infty)$.
It follows that the integral term in (\TEneg) is bounded by $CT^{-n/2k}$, hence converges to $0$ as $T\to\infty$.
Consequently, $E(w_{0,T}(-\log T))<0$ for $T$ sufficiently large.
By Proposition~4.3, this contradicts the global existence of $u$
and we conclude that $K=0$.

Finally, in view of (\TEnegA) or (\TEnegB) with $K=0$, Theorem~1 is now a consequence of Corollary~3. \fin
\smallskip

\smallskip
{\bf Remark~7.1.}  Let us justify the assertion in Remark~3.1(iii).
From the proof of Corollary~3, under assumption (\SdecayHaa), we actually get
$${\cal N}(u_0,t)\le Ct^{-\eta},\quad t\ge 1,$$
with some $\eta=\eta(n,p,q)>0$, where ${\cal N}(u_0,t)$ is defined in Proposition~5.1.
By Proposition~5.1(ii), it follows
that $\|u(t)\|_{M^{2,4/(p-1)}}\le Ct^{-\eta/(p+1)}$. For $t\ge t_1$ large enough, we deduce from Proposition~6.1 that
$$\|u(t)\|_\infty\le C_0 \, (t/2)^{-1/(p-1)}\|u(t/2)\|_{M^{2,4/(p-1)}}\le Ct^{-\beta},\quad t\ge t_1,$$
where $\beta={1\over p-1}+{\eta\over p+1}>{1\over p-1}$.
Consequently $h(t):=\int_0^t\|u(s)\|_\infty^{p-1}\,ds$ satisfies $\displaystyle\sup_{t>0}h(t)<\infty$ 
and $z(t):=e^{h(t)}e^{-tA}|u_0|$ solves
$$z_t-\Delta z=h'(t)z=\|u(t)\|_\infty^{p-1}z \ge |u|^{p-1}z.$$
Therefore, $u\le z$ by the maximum principle, and similarly $u\ge -z$.
The claimed property (\LinearBehavior) follows.
On the other hand, by the proof of Theorem~1, under assumption (\SdecayHaaa), 
we may assume $u_0\in L^{q^*}(\Rn)$ after a time shift.
We thus again obtain (\LinearBehavior) and then (\LinearBehaviorB).
\fin
\smallskip

\medskip

{\bf 8. Proof of Theorem 4.}

\medskip

It is based on the following continuous dependence property 
with respect to the critical Morrey norms.

\proclaim Proposition~8.1. 
Let $\Omega\subset\Rn$ be a (possibly unbounded) domain of class $C^{2+\alpha}$.
Let $q\in (1,\infty)$ and $p\ge 1+{2q\over n}$.
For any $u_0\in L^\infty\cap M^{q,2q/(p-1)}(\Omega)$, and any $T_0\in (0,T_{max}(u_0))$,
there exist $\delta, M>0$ (depending on $u_0$ and $T_0$) with the following property.
If $v_0\in L^\infty\cap M^{q,2q/(p-1)}(\Omega)$ satisfies
$$\|u_0-v_0\|_{M^{q,2q/(p-1)}}\le \delta,$$
then $T_{max}(v_0)>T_0$ and the corresponding solution $v$ of (\equp) satisfies
$$\|u(t)-v(t)\|_{M^{q,2q/(p-1)}}\le M \|u_0-v_0\|_{M^{q,2q/(p-1)}},\quad 0\le t\le T_0.
\eqno(\depcontMorrey)$$

\medskip

{\bf Proof.} It relies on suitable modifications of the proof of Proposition~6.1.
We set $\lambda=2q/(p-1)$ and denote again $u^p=|u|^{p-1}u$ and $|\cdot |_m=\|\cdot\|_{M^{m,\lambda}(\Omega)}$. 
In all the proof, $C_i$ (resp., $C$) denote fixed (resp., generic) positive constants depending only on $n,p,q,r$.
\smallskip

{\bf Step 1.} Fix $r$ such that $1<r/p<q<r$ and let $\beta={1\over p-1}-{\lambda\over 2r}>0$. 
Let $C_1$ be the constant from (\SMorreyInvEstA). 
We claim that there exists $\delta_0=\delta_0(n,p,q,r)>0$ with the following property:
for any $\delta\in (0,\delta_0]$, there exists $t_1=t_1(\delta)\in (0,T_0)$ (depending on $u_0$) such that
$$
|u_0-v_0|_q\le \delta\ \ \Longrightarrow\ \  
\left\{\eqalign{
&t^\beta |u(t)-v(t)|_r \le 2C_1|u_0-v_0|_q\cr
&\hbox{for all $t\in \bigl(0,\min(t_1(\delta),T_{max}(v_0))\bigr)$}.}\right.
\eqno(\supsbetaA)$$

To prove (\supsbetaA), let $\delta>0$ and assume $|u_0-v_0|_q\le \delta$. Since $u_0\in M^{q,\lambda}\cap L^\infty\subset M^{r,\lambda}$, inequalities (\ContractionMorrey) and (\SMorreyInvEstAB) guarantee the existence of 
$t_1=t_1(\delta)\in (0,T_0)$ such that  
$$s^\beta |u(s)|_r+s^{1/(p-1)} \|u(s)\|_\infty \le \delta, \qquad 0<s\le t_1(\delta).
\eqno(\supsbeta)$$
Set $\tilde t_1=\tilde t_1(\delta,v_0)=\min(t_1(\delta),T_{max}(v_0))$ and let 
$$\tau=\sup\,\bigl\{t<\tilde t_1;\, s^\beta |u(s)-v(s)|_r\le 2C_1 |u_0-v_0|_q\ \hbox{ for all $s\in(0,t)$}\bigr\}$$
(note that $\tau>0$, owing to (\SMorreyInvEstAB) and (\ContractionMorrey) applied to $u_0,v_0\in M^{r,\lambda}$
with some $t_0<\min(T_{max}(u_0),$ $T_{max}(v_0))$).
For all $s\in (0,\tau)$, we have $s^\beta |v(s)|_r\le (2C_1+1)\delta$ hence,
using (\PropMorreyA), (\PropMorreyB),
$$\eqalign{ |u^p(s)-v^p(s)|_{r/p}
&\le p\Bigl| \bigl(|u(s)|^{p-1}+|v(s)|^{p-1}\bigr) |u(s)-v(s)|\Bigr|_{r/p} \cr
&\le p\bigl(|u(s)|_r^{p-1}+|v(s)|_r^{p-1}\bigr)|u(s)-v(s)|_r \cr
\noalign{\vskip 1mm}
&\le C\delta^{p-1}|u_0-v_0|_qs^{-\beta p}.}
\eqno(\upvprp)$$
Therefore, by (\SMorreyInvEstA), (\SMorreyInvEstAA) and (\intinv), we deduce that, for all $t\in (0,\tau)$,
$$\eqalign{
t^\beta|u(t)-v(t)|_r
&\le t^\beta|e^{-tA}(u_0-v_0)|_r+t^\beta\int_0^t |e^{-(t-s)A}(u^p(s)-v^p(s))|_r\, ds \cr
&\le C_1|u_0-v_0|_q+C_1t^\beta\int_0^t (t-s)^{-\lambda(p-1)/2r} |u^p(s)-v^p(s)|_{r/p}\, ds \cr
&\le C_1|u_0-v_0|_q+C\delta^{p-1}|u_0-v_0|_qt^\beta\int_0^t (t-s)^{-q/r}s^{-\beta p}\, ds \cr
&\le (C_1+C\delta^{p-1})|u_0-v_0|_q. }$$
Assuming $\delta\in(0,\delta_0]$ with $\delta_0=\delta_0(n,p,q,r)>0$ 
small enough and arguing as before (\SMorreyInvEst), we deduce that $\tau=\tilde t_1$,
hence the claim.

\smallskip
{\bf Step 2.} We claim that there exist $\delta_1=\delta_1(n,p,q,r)\in (0,\delta_0]$ and $C_3, C_4>0$ such that, 
if $|u_0-v_0|_q\le \delta_1$, then
$$T_{max}(v_0)>t_1(\delta_1),\qquad
t^{1/(p-1)}\|u(t)-v(t)\|_\infty\le C_3 |u_0-v_0|_q,\quad t\in (0,t_1(\delta_1)],
\eqno(\supsbetaB)$$
where $t_1(\delta_1)$ is given by Step~1, and
$$|u(t)-v(t)|_q\le C_4 |u_0-v_0|_q,\quad t\in (0,t_1(\delta_1)].
\eqno(\supsbetaBB)$$

To prove (\supsbetaB), for $C_3>0$ and $\delta_1\in (0,\delta_0]$ to be chosen later, 
we assume $|u_0-v_0|_q\le \delta_1$. Denote
 $\tilde t_1=\tilde t_1(\delta_1,v_0)=\min(t_1(\delta_1),T_{max}(v_0))$ and let 
$$\tau_1=\sup\,\bigl\{t<\tilde t_1;\, s^{1/(p-1)} \|u(s)-v(s)\|_\infty\le C_3|u_0-v_0|_q \ \hbox{for all $s\in(0,t)$}\bigr\}$$
(note that $\tau_1>0$, since $\|u(s)\|_\infty$ and $\|v(s)\|_\infty$ are bounded on $[0,t]$ for some small $t>0$).
In particular, by (\supsbeta), we have $s^{1/(p-1)} (\|u(s)\|_\infty\vee\|v(s)\|_\infty)\le (C_3+1)\delta_1$ on $(0,\tau_1)$.
For all $s\in(0,\tau_1)$, we have
$$\eqalign{
\|u^p(s)-v^p(s)\|_\infty
&\le p\bigl(\|u(s)\|_\infty^{p-1}+\|v(s)\|_\infty^{p-1}\bigr)\|u(s)-v(s)\|_\infty \cr
&\le 2p(C_3+1)^p\delta_1^{p-1}|u_0-v_0|_q s^{-{p\over p-1}}.}$$
Using (\Varconsthalf), (\SMorreyInvEstB), (\supsbetaA) and recalling $\beta={1\over p-1}-{\lambda\over 2r}$, 
we deduce that, for all $t\in(0,\tau_1)$,
$$\eqalign{  
t^{1\over p-1}\|u(t)-v(t)\|_\infty 
&\leq t^{1\over p-1} \bigl\|e^{-{t\over 2}A} (u(\textstyle{t\over 2})-v(\textstyle{t\over 2})) \bigr\|_\infty
+t^{1\over p-1}\displaystyle\int_{t/2}^t \bigl\|u^p(s)-v^p(s)\|_\infty\, ds\cr
&\leq 2^{\lambda\over 2r}C_2t^{{1\over p-1}-{\lambda\over 2r}} |u(\textstyle{t\over 2})-v(\textstyle{t\over 2})|_r
+2^{p\over p-1}p(C_3+1)^p\delta_1^{p-1}|u_0-v_0|_q \cr
\noalign{\vskip 1mm}
& \le  \Bigl[2^{{\lambda\over 2r}+\beta+1}C_1C_2+2^{p\over p-1}p(C_3+1)^p\delta_1^{p-1}\Bigr] |u_0-v_0|_q.
}$$
Choosing $C_3=2^{{\lambda\over 2r}+\beta+2}C_1C_2$ and $\delta_1=\delta_1(n,p,q,r)>0$ small enough, 
we get $\tau_1=\tilde t_1$, hence (\supsbetaB).

To prove (\supsbetaBB), we set $\gamma={\lambda\over 2}\bigl({p\over r}-{1\over q}\bigr)>0$.
Using Proposition~4.1 and (\upvprp), noting that $\gamma+\beta p=1$, we obtain
$$\eqalign{  
|u(t)-v(t)|_q
&\leq C|e^{-tA} (u_0-v_0)|_q
+C\displaystyle\int_0^t (t-s)^{-\gamma}|u^p(s)-v^p(s)|_{r/p}\, ds\cr
&\leq C|u_0-v_0|_q
+C\delta_1^{p-1}|u_0-v_0|_q\displaystyle\int_0^t (t-s)^{-\gamma}s^{-\beta p}\, ds
\leq C_4|u_0-v_0|_q.}$$

{\bf Step 3.} Let $M_0=\sup_{t\in [0,T_0]}\|u(t)\|_\infty<\infty$, 
let $T_1:=t_1(\delta_1)$ be given by Steps~1 and 2, and assume 
$$|u_0-v_0|_q<\min\bigl(\delta_1,C_3^{-1}T_1^{1/(p-1)}\bigr).$$
Then, by (\supsbetaB), we have $T_2:=\min(T_0,T_{max}(v_0))>T_1$ and
$$\tau_2:=\sup\,\bigl\{t\in (T_1,T_2);\, \|u(s)-v(s)\|_\infty\le 1\ \hbox{for all $s\in [T_1,t]$}\bigr\}>T_1.$$
Next, by the maximum principle, 
using $\max(\|u(s)\|_\infty,\|v(s)\|_\infty) \le M_0+1$ on $[T_1,\tau_2)$, we easily~get
$$|u(t)-v(t)|\le e^{p(M_0+1)^{p-1}(t-T_1)} \bigl|e^{-(t-T_1)A}(u(T_1)-v(T_1))\bigr| 
\quad\hbox{ in $\Omega\times[T_1,\tau_2)$.}
\eqno(\supsbetaC)$$
Setting $M_1=e^{p(M_0+1)^{p-1}(T_0-T_1)}$, $M_2=M_1C_3T_1^{-1/(p-1)}$
and using (\supsbetaB), we in particular obtain
$$\|u(t)-v(t)\|_\infty\le M_1\|u(T_1)-v(T_1)\|_\infty 
\le M_2|u_0-v_0|_q \quad\hbox{ in $\Omega\times[T_1,\tau_2)$.}$$
Now further assuming $|u_0-v_0|_q<\delta:=\min\bigl(\delta_1,C_3^{-1}T_1^{1/(p-1)},(2M_2)^{-1}\bigr)$, 
it follows that $\tau_2=T_2$.
Consequently, $T_{max}(v_0)>T_0$. 
Finally, going back to (\supsbetaC) and using (\ContractionMorrey) and (\supsbetaBB), we conclude that
$$|u(t)-v(t)|_q\le M_1 |e^{-(t-T_1)A}(u(T_1)-v(T_1))|_q
\le M_1 |u(T_1)-v(T_1)|_q\le C_4M_1 |u_0-v_0|_q,$$
for all $t\in [T_1,T_0]$. This combined with (\supsbetaBB) proves (\depcontMorrey).
\fin

\medskip

{\bf Proof of Theorem 4.}
(i) Let $\eps_0$ be given by Proposition~6.1 with $q=2$.
By assumption (\hypdecayuzeto) and Proposition~5.1, 
since $u$ is assumed to be global, there exists $T_0>0$ such that 
$|u(T_0)|_2\le \eps_0/2$, where $|\cdot |_2=\|\cdot\|_{M^{2,4/(p-1)}}$.
\smallskip

By Proposition~4.2, we have $u_0\in M^{2,4/(p-1)}\cap L^\infty$.
By Proposition~8.1, there exists $\eta>0$ such that if $|u_0-v_0|_2\le \eta$,
then $T_{max}(v_0)>T_0$ and $|u(T_0)-v(T_0)|_2\le\eps_0/2$, hence $|v(T_0)|_2\le\eps_0$.
We then conclude from Proposition~6.1 that $v$ is global, which proves (\GopenA).
The inclusions (\GopenB) and (\GopenC) then follow from Proposition~4.2.

\smallskip
(ii) Assume for contradiction that $T_{max}(\lambda^*\phi)=\infty$.
Then, since $|\phi|_2<\infty$ by (\hypdecayuzeto) and (\compMorreyGnorm), we would have $T_{max}(\lambda\phi)=\infty$
for $\lambda\to\lambda^*_+$ by assertion (i), contradicting the definition of $\lambda^*$. \fin

\medskip

\goodbreak
\centerline{\bf Appendix.}
\medskip

In this appendix, for convenience and self-containedness, we provide proofs of some important known results 
that we have used,
namely Proposition~4.1 from [\Kato] and Proposition~4.3 from [\GKa, \GKb].
We also give the proofs of Lemmas~7.1 and 7.2.

\medskip
{\bf Proof of Proposition~4.1.}
By the maximum principle, it suffices to prove the assertion when $\Omega=\R^n$. 
Set $v(\cdot,t)=e^{-tA}f=G_t*f$. 
Let $p\in [1,\infty)$ (the case $p=q=\infty$ is obvious due to $M^{\infty,\lambda}=L^\infty$). 
By Jensen's inequality and $\|G_t\|_1=1$, we have
$|v(\cdot,t)|^p\le G_t*|f|^p$. By Fubini's theorem, for any $R>0$, it follows that
$$\eqalign{
\int_{B_R}|v(x,t)|^p\,dx
&\le \int_{B_R}\Bigl(\int_{\R^n}G_t(y)|f(x-y)|^p\, dy\Bigr)dx \cr
&=\int_{\R^n}G_t(y) \Bigl(\int_{B_R}|f(x-y)|^p\,dx\Bigr) dy
=\int_{B_R(x)}|f|^p\,dx \le R^{n-\lambda}\|f\|^p_{M^{p,\lambda}}.}$$
Consequently, $\|e^{-tA}f\|^p_{M^{p,\lambda}}\le \|f\|^p_{M^{p,\lambda}}$, 
hence (\ContractionMorrey).
\smallskip

Next set $\rho(r)=\int_{B_r}|f|^p\,dx$.
We claim that, for each $R>0$ and $\phi\in C^1([0,R])$,
$$\int_{B_R}\phi(|x|)|f(x)|^p\,dx=\phi(R)\rho(R)-\int_0^R \phi'(r)\rho(r)\,dr.
\eqno(\ClaimKatoMorrey)
$$
By density of $C(\overline B_R)$ in $L^p(B_R)$, 
it suffices to prove this when $f\in C(\overline B_R)$. Set $h(r)=\int_{\partial B_r}|f|^p\,d\sigma_r$,
where $d\sigma_r$ is the surface measure on $\partial B_r$.
Then $h\in C([0,R])$ and $\rho(r)=\int_0^r h(s) ds$. Integrating by parts, we get
$$\eqalign{
\int_{B_R}\phi(|x|)|f(x)|^p\,dx
&=\int_0^R \phi(r)h(r)\,dr =\int_0^R \phi(r)\rho'(r)\,dr\cr
&=\phi(R)\rho(R)-\int_0^R \phi'(r)\rho(r)\,dr.}$$
For fixed $t>0$, using H\"older's inequality and then 
applying (\ClaimKatoMorrey) with $\phi(r)=K_t(r)=(4\pi t)^{-n/2}e^{-r^2/4t}$, we obtain
$$\eqalign{|v(0,t)|^p
&\le \Bigl(\int_{B_R} K_t(|x|)|f(x)|\,dx\Bigr)^p \le\int_{B_R} K_t(|x|)|f(x)|^p\,dx\cr
&\le K_t(R)\rho(R)-\int_0^R K_t'(r)\rho(r)\,dr\cr
&\le K_t(R)R^{n-\lambda}\|f\|^p_{M^{p,\lambda}}
+(4\pi t)^{-n/2}(2t)^{-1}\|f\|^p_{M^{p,\lambda}}\int_0^R r^{1+n-\lambda}e^{-r^2/4t}\,dr.}$$
Letting $R\to\infty$, we get
$$|v(0,t)|^p\le t^{-1-(n/2)}\|f\|^p_{M^{p,\lambda}}\int_0^\infty r^{1+n-\lambda}e^{-r^2/4t}\,dr 
\le C(n,\lambda)t^{-\lambda/2}\|f\|^p_{M^{p,\lambda}}.$$
By translation invariance, the same holds at any point $x_0\in \R^n$ instead of $0$.
This yields~(\RegulMorrey) for $q=\infty$. 
The general case follows by interpolating between the cases $q=p$ and $q=\infty$.
\fin
\medskip

{\bf Proof of Proposition~4.3.}
Assume $a=0$ without loss of generality.
For all $s>s_0$, we compute
$$
{d\over ds}\int_{D(s)}|w|^q\rho=\int_{D(s)}\rho\,\partial_s(|w|^q)
+e^{s/2}\int_{\partial D(s)}|w|^q\rho {y\cdot\nu\over |y|}\,d\sigma,
\eqno(\SGKdds)
$$
where $q=2$ or $p+1$, and the boundary term vanishes since $w=0$ on $\partial D(s)$.
Here and below, for $s>s_0$, all the procedures are justified owing to the fast decay of the Gaussian weight $\rho$
and parabolic regularity. By integration by parts, we have
$$\eqalign{
  {1\over 2}\,{d\over ds}\int_{D(s)}w^2\rho
    &=\int_{D(s)}ww_s\rho
    =\int_{D(s)}w\bigl[\nabla\cdot(\rho\nabla w)+\rho |w|^{p-1}w-\beta\rho w\bigr] \cr
    &=\int_{D(s)}\bigl[-|\nabla w|^2-\beta w^2+|w|^{p+1}\bigr]\rho
    = -2E(w)+{p-1\over p+1}\int_{D(s)}|w|^{p+1}\rho}
$$
i.e., (\QGKidA). Next, using
$$
{d\over ds}\int_{D(s)}|\nabla w|^2\rho
=-{d\over ds}\int_{D(s)}w\nabla\cdot(\rho\nabla w)
 $$
and noting that the variation of the domain again does not produce a boundary term, due to $w=0$ on $\partial D(s)$,
we obtain
$$\eqalign{
{d\over ds}\int_{D(s)}|\nabla w|^2\rho
&=-\int_{D(s)}w_s\nabla\cdot(\rho\nabla w)-\int_{D(s)}w\nabla\cdot(\rho\nabla w_s) \cr
&=-\int_{D(s)}w_s\nabla\cdot(\rho\nabla w)+\int_{D(s)}\rho(\nabla w\cdot\nabla w_s) \cr
&=-2\int_{D(s)}w_s\nabla\cdot(\rho\nabla w)+\int_{\partial D(s)}\rho w_s w_\nu\,d\sigma. 
}
\eqno(\SGKconvex)
$$
Since 
$$w_s=-\beta w-{y\over 2}\cdot\nabla w+e^{-(\beta+1)s}u_t(a+e^{-s/2}y,T-e^{-s})$$
and since the tangential derivatives of $w$ vanish on $\partial D(s)$,
we have $w_s=-{y\over 2}\cdot\nabla w=-{1\over 2}(y\cdot\nu)w_\nu$ on $\partial D(s)$.
Therefore, the integrand of the boundary term in (\SGKconvex) can be written as
$-{1\over 2}\rho (y\cdot\nu)|w_\nu|^2\le 0$, owing to the starshapedness of $\Omega$ with respect to~$a$. 
Consequently,
$${1\over 2}{d\over ds}\int_{D(s)}|\nabla w|^2\rho
\le -\int_{D(s)}w_s\nabla\cdot(\rho\nabla w). $$
On the other hand, by (\SGKdds),  we have
$$ {d\over ds}\int_{D(s)}
\Bigl({\beta\over 2} w^2-{1\over p+1}|w|^{p+1}\Bigr)\rho
 = \int_{D(s)}(\beta w-|w|^{p-1}w)w_s\rho.$$
Summing the last two formulas and using equation (\QGKeqrho), 
we obtain (\QGKidB).
Since the continuity of $E(s)$ at $s=s_0$ is guaranteed by the assumption $u_0\in BC^1(\Omega)$,
we deduce in particular that $E(s)\le E(s_0)$ for all $s>s_0$.

Next denote $\psi(s):=\int_{D(s)}w^2(s)\rho$.
Then (\QGKidA), Jensen's inequality and (\QGKidB) imply
$$ {1\over 2}{d\psi\over ds}\geq-2E(w(s))+C(n,p)\psi^{(p+1)/2}(s)
\geq -2E(w(s_0))+C(n,p)\psi^{(p+1)/2}(s).
$$
This guarantees (\QGKenergypos) and (\QGKLtwocontrol) 
(otherwise $\psi$ has to blow up in finite time).
\fin

   \def\Xint#1{\mathchoice
      {\XXint\displaystyle\textstyle{#1}}%
      {\XXint\textstyle\scriptstyle{#1}}%
      {\XXint\scriptstyle\scriptscriptstyle{#1}}%
      {\XXint\scriptscriptstyle\scriptscriptstyle{#1}}%
      \!\int}
   \def\XXint#1#2#3{{\setbox0=\hbox{$#1{#2#3}{\int}$}
        \vcenter{\hbox{$#2#3$}}\kern-.5\wd0}}
   
   \def\dashint{\Xint-}

\medskip
{\bf Proof of Lemma~7.1.}
(i) Fix a cut-off function $\rho\in C^\infty_0([0,\infty))$, with $0\le \rho\le 1$,
 such that $\rho(s)=1$ for $s\le 1$ and $\rho(s)=0$ for $s\ge 2$.
For each integer $j\ge 1$, set $A_j=B_{2j}\setminus \overline B_j$ and define
$$v_j(x)=\rho\bigl(\textstyle{|x|\over j}\bigr)\Bigl(v(x)-\displaystyle\dashint_{A_j}v\Bigr).
\eqno(\vjrho)$$
Since $\rho$ is compactly supported, it is clear that $v_j\in W^{1,q}(\Rn)$. We compute
$$\nabla v_j(x)-\nabla v(x)={1\over j}{x\over |x|}\rho'\bigl(\textstyle{|x|\over j}\bigr)\Bigl(v(x)-\displaystyle\dashint_{A_j}v\Bigr)
+\Bigl(\rho\bigl(\textstyle{|x|\over j}\bigr)-1\Bigr)\nabla v(x),$$
hence
$$\|\nabla v_j-\nabla v\|_{L^q(\Rn)}\le {C\over j}\Bigl\|v-\displaystyle\dashint_{A_j}v\Bigr\|_{L^q(A_j)}
+\|\nabla v\|_{L^q(\Rn\setminus B_j)}.$$

By the Poincar\'e-Wirtinger inequality, there exists a constant $C_0=C_0(n,q)>0$ such that
$$\Bigl\|\phi-\displaystyle\dashint_{A_1}\phi\Bigl\|_{L^q(A_1)}\le 
C_0\|\nabla\phi\|_{L^q(A_1)}\quad\hbox{ for all $\phi\in W^{1,p}(A_1)$}.$$
By a simple scaling argument, it follows that
$$\Bigl\|\phi \,-\, \displaystyle\dashint_{A_j}\phi\Bigl\|_{L^q(A_j)}\le 
C_0 {\hskip 1pt} j \|\nabla\phi\|_{L^q(A_j)}\quad\hbox{ for all $\phi\in W^{1,p}(A_j)$,\ \ $j\ge 1$}.$$
Therefore,
$$\|\nabla v_j-\nabla v\|_{L^q(\Rn)}\le CC_0 \|\nabla v\|_{L^q(A_j)}
+\|\nabla v\|_{L^q(\Rn\setminus B_j)}
\le (1+CC_0) \|\nabla v\|_{L^q(\Rn\setminus B_j)}.$$
Since the RHS goes to $0$ as $j\to\infty$, it follows that $\nabla v_j\to \nabla v$ in $L^q(\Rn)$.

Next, let $K_j:=\dashint_{A_j}v$. Since $|K_j|\le \|v\|_\infty$, we may assume,
after extracting a subsequence, that $\lim_{j\to\infty} K_j=K\in \R$.
Also, since $v_j\in W^{1,q}(\Rn)$, we have 
$$\|v_j\|_{L^{q^*}(\Rn)} \le S\|\nabla v_j\|_{L^q(\Rn)},
\eqno(\vjSobolev)$$
where $S=S(n,q)>0$ is the Sobolev constant.
Observe that, by (\vjrho), for all given $R>0$, $x\in B_R$ and all $j\ge R$, we have $v_j(x)+K_j=v(x)$,
hence
$$v(x)-K=v_j(x)+(K_j-K),\quad x\in B_R,\ j\ge R.$$
Therefore, by (\vjSobolev),
$$
\eqalign{
\|v-K\|_{L^{q^*}(B_R)} 
&\le \|v_j\|_{L^{q^*}(B_R)}+|B_R|^{1/q^*}|K_j-K|  \cr
&\le S\|\nabla v_j\|_{L^q(\Rn)}+|B_R|^{1/q^*}|K_j-K|.}$$
Passing to the limit $j\to\infty$, we obtain
$$\|v-K\|_{L^{q^*}(B_R)} \le S\|\nabla v\|_{L^q(\Rn)}.$$
Uppon letting $R\to\infty$, we conclude that $v-K\in L^{q^*}(\Rn)$.

\medskip

(ii) Let $\eps>0$. By assumption, there exists $A>0$ such that
$$|\nabla v(x)|\le \beta\eps |x|^{-\beta-1},\quad |x|\ge A.
\eqno(\DecayGrad)$$
For each $\omega\in S^{n-1}$ and $s>r\ge A$, we compute
$$
\eqalign{|v(s\omega)-v(r\omega)|
&=\Bigl|\int_0^1 (s-r)\omega\cdot\nabla v\bigl((r+t(s-r))\omega\bigr)\,dt\Bigr| \cr
&\le \int_0^1 (s-r)\bigl|\nabla\phi\bigl((r+t(s-r))\omega\bigr)\bigr|\,dt\cr
&\le\beta\eps \int_0^1 (s-r) (r+t(s-r))^{-\beta-1}\,dt
\le\eps r^{-\beta}.}$$
By the Cauchy property, we deduce that $\ell(\omega)=\lim_{r\to\infty} v(r\omega)\in\R$ exists.
Moreover, by letting $s\to\infty$ in the previous inequality, we obtain
$$
|v(r\omega)-\ell(\omega)|\le\eps r^{-\beta},\quad r\ge A.$$
To prove the assertion, it thus suffices to show that $\ell(\omega)$ is independent of $\omega$.
Consider $\omega_1,\omega_2\in S^{n-1}$ such that $\langle\omega_1,\omega_2\rangle\ge 0$.
We see that, for all $r>0$ and $t\in [0,1]$,
$$\eqalign{
|r\omega_1+rt(\omega_2-\omega_1)|^2
&=r^2\bigl[1+2t^2(1-\langle\omega_1,\omega_2\rangle)+2t(\langle\omega_1,\omega_2\rangle-1)\bigr]  \cr
&=r^2\bigl[1-2t(1-t)+2t(1-t)\langle\omega_1,\omega_2\rangle\bigr] \ge {1\over 2}\,r^2.
}$$
For $r\ge A\sqrt{2}$, it follows from (\DecayGrad) that
$$
\eqalign{|v(r\omega_2)-v(r\omega_1)|
&=\Bigl|\int_0^1 r(\omega_2-\omega_1)\cdot\nabla v\bigl(r\omega_1+rt(\omega_2-\omega_1)\bigr)\,dt\Bigr| \cr
&\le 2r\int_0^1 \bigl|\nabla v\bigl(r\omega_1+rt(\omega_2-\omega_1)\bigr)\bigr|\,dt
\le 2\beta \eps r\Bigl({r\over \sqrt 2}\Bigr)^{-\beta-1}=\beta\eps 2^{{\beta+3\over 2}}r^{-\beta},}$$
hence $\ell(\omega_1)=\ell(\omega_2)$ by letting $r\to\infty$.
Denoting by $(e_i)_{1\le i\le n}$ the canonical basis of $\Rn$, 
we have in particular $\ell(e_1)=\ell(e_2)=\ell(-e_1)$. For any $\omega\in S^{n-1}$, since 
either $\langle\omega,e_1\rangle\ge 0$ or $\langle\omega,-e_1\rangle\ge 0$,
it follows that $\ell(\omega)=\ell(e_1)$ and the assertion is proved.
 \fin

\medskip
{\bf Proof of Lemma~7.2.}
Fix $\tau\in (0,\min(1,T_{max}))$ and 
let $M=\sup_{t\in (0,\tau)} \|u(t)\|_\infty$. Using
$$\nabla u(t)=\nabla G_t*u_0+\int_0^t \nabla G_{t-s}*u^p(s)\, ds, 
\eqno(\varconstNabla)$$
we first obtain, for all $t\in (0,\tau)$,
$$\|\nabla u(t)\|_\infty\le Ct^{-1/2}\|u_0\|_\infty+CM^p\int_0^t (t-s)^{-1/2}\, ds
\le C_1t^{-1/2}.
\eqno(\varconstNablaB)$$
Next rewriting (\varconstNabla) as
$\nabla u(t)=G_t*\nabla u_0+\int_0^t G_{t-s}*\nabla u^p(s)\, ds$,
we get
$$|\nabla u(t)|\le G_t*|\nabla u_0|+N\int_0^t G_{t-s}*|\nabla u(s)|\, ds, $$
where $N=pM^{p-1}$.
Next set $z(t):=e^{Nt}G_t*|\nabla u(t)|$. Noting that $\|\nabla u(t)\|_\infty\in L^1(0,\tau)$ owing to
(\varconstNablaB), we have, by direct computation,
$$z(t)=G_t*|\nabla u_0|+N\int_0^t G_{t-s}*z(s)\, ds.$$
Therefore,
$$\bigl(|\nabla u(t)|-z(t)\bigr)_+\le N\int_0^t G_{t-s}*\bigl(|\nabla u(s)|-z(s)\bigr)_+\, ds,
\quad 0<t<\tau,$$
hence
$$\bigl\| \bigl(|\nabla u(t)|-z(t)\bigr)_+\bigr\|_\infty \le N\int_0^t \bigl\|\bigl( |\nabla u(s)|-z(s)\bigr)_+\bigr\|_\infty\, ds,
\quad 0<t<\tau.$$
Since $\|(|\nabla u(t)|-z(t))_+\|_\infty \in L^1(0,\tau)$, we may apply Gronwall's lemma
to conclude that $\|(|\nabla u(t)|-z(t))_+\|_\infty=0$, hence $|\nabla u(t)|\le z(t)$ on $(0,\tau)$,
 which yields the desired conclusion
for possibly smaller $\tau$.
\fin

\vskip 2mm
{\baselineskip=11pt \parindent=0.7cm

\font\rmn=cmr9
\font\sln=cmsl9
\font\rmb=cmbx8 scaled 1125 \rm

\rmn \eightpoint

\bigskip
\centerline{\bf REFERENCES}
\medskip\medskip

\item{[\BKZ]} 
P. Biler, G. Karch and J. Zienkiewicz,
Optimal criteria for blowup of radial and $n$-symmetric solutions of chemotaxis systems,
{\sln Nonlinearity} 28 (2015), 4369.

\smallskip

\item{[\BS]} 
S. Blatt and M. Struwe,
An analytic framework for the supercritical Lane-Emden equation and its gradient flow,
{\sln Int. Math. Res. Notices}
2015 (2015), 2342-2385.
\smallskip

\item{[\CFG]}
X. Chen, M. Fila and J.-S. Guo,
Boundedness of global solutions of a supercritical parabolic equation,
{\sln Nonlinear Anal.} 68 (2008), 621-628.
\smallskip

\item{[\CDZ]}
K.-S. Chou, S.-Z. Du, and G.-F. Zheng.
On partial regularity of the borderline solution of semilinear parabolic problems,
{\sln Calc. Var. Partial Differential Equations} 30 (2007), 251-275.
\smallskip

\item{[\Fi]} 
M. Fila,
Boundedness of global solutions of nonlinear diffusion equations,
{\sln J. Differential Equations}
98 (1992), 226-240.

\smallskip
\item{[\FK]}
M. Fila and J.R. King,
Grow up and slow decay in the critical Sobolev case,
{\sln Netw. Heterog. Media}  7 (2012), 661-671.

\smallskip
\item{[\FM]}
M. Fila and N. Mizoguchi,
Multiple continuation beyond blow-up,
{\sln Diff. Int. Equations} 20 (2007), 671-680.

\smallskip
\item{[\GK]}
 V.A. Galaktionov and J.R. King,
Composite structure of global unbounded solutions of nonlinear heat equations with critical Sobolev exponents,
{\sln J. Differential Equations} 189 (2003), 199-233.
\smallskip

\item{[\GV]}
 V.A. Galaktionov and J.L. V\'azquez,
Continuation of blow-up solutions of nonlinear heat equations in several space dimensions,
{\sln Comm. Pure Appl. Math.} 50 (1997), 1-67.
\smallskip

\item{[\Gig]}
 Y. Giga,
A  bound for global solutions of semilinear heat equations,
{\sln Comm. Math. Phys.}
103 (1986), 415-421.
\smallskip

\item{[\GKa]}
Y. Giga and R.V. Kohn,
Asymptotically self-similar blow-up of semilinear heat equations,
{\sln Comm. Pure Appl. Math.} 38 (1985), 297-319.
\smallskip

\item{[\GKb]}
Y. Giga and R.V. Kohn,
Characterizing blowup using similarity variables,
{\sln Indiana Univ. Math. J.} 36 (1987), 1-40.
\smallskip

\item{[\GMi]}
Y. Giga and T. Miyakawa,
Navier-Stokes flow in $\R^3$ with measures as initial vorticity and Morrey spaces,
{\sln Comm. Partial Differential Equations} 14 (1989), 577-618. 

\smallskip

\item{[\HW]}
A. Haraux and F.B. Weissler,
Non-uniqueness for a semilinear initial value problem,
{\sln Indiana Univ. Math. J.} 31 (1982), 167-189.
\smallskip

\item{[\Kar]}
G. Karch,
Scaling in nonlinear parabolic equations,
{\sln J. Math. Anal. Appl.} 234 (1999), 534-558.

\smallskip

\item{[\Kato]}
T. Kato, 
Strong Solutions of the Navier-Stokes Equation in Morrey Spaces,
{\sln Bol. Soc. Bras. Mat.}
22 (1992), 127-155.
\smallskip

\item{[\Kav]}
O. Kavian, 
Remarks on the large time behaviour of a nonlinear diffusion equation, 
{\sln Ann. Inst. H. Poincar\'e Analyse Non Lin\'eaire} 4 (1987), 423-452.

\smallskip

\item{[\Kaw]}
T. Kawanago,
Asymptotic behavior of solutions of a semilinear heat equation with subcritical nonlinearity,
{\sln Ann. Inst. H. Poincar\'e Anal. non lin\'eaire} 13 (1996), 1-15.
\smallskip

\item{[\KS]}
H. Kozono and H. Sohr, 
New a priori estimates for the Stokes equations in exterior domains,
{\sln Indiana Univ. Math. J.} 40 (1991), 1-27.
\smallskip

\item{[\LN]}
 T. Lee and W.-M. Ni,
Global existence, large time
behaviour and life span of solutions of a semilinear parabolic Cauchy problem,
{\sln Trans. Amer. Math. Soc.} 333 (1992), 365-378.
\smallskip

\item{[\LR]}
P.-G. Lemari\'e-Rieusset,
Small data in an optimal Banach space for the parabolic-para\-bolic and parabolic-elliptic 
Keller-Segel equations in the whole space,
{\sln Adv. Differential Equations} 18 (2013), 1189-1208.
\smallskip

\item{[\MMCPAM]}
H. Matano and F. Merle,
On nonexistence of type II blowup for a supercritical nonlinear heat equation, 
{\sln Comm. Pure Appl. Math.} 57 (2004), 1494-1541.
\smallskip

\item{[\MMJFA]}
H. Matano and F. Merle,
Classification of type I and type II behaviors for a supercritical nonlinear heat equation,
{\sln J. Funct. Anal.} 256 (2009), 992-1064.

\smallskip

\item{[\MMJFAb]}
H. Matano and F. Merle,
Threshold and generic type I behaviors for a supercritical nonlinear heat equation,
{\sln J. Funct. Anal.} 261 (2011) 716-748.

\smallskip

\item{[\Miz]}
 N. Mizoguchi,
On the behavior of solutions for a semilinear parabolic equation with supercritical nonlinearity,
{\sln Math. Z.} 239 (2002), 215-229.
\smallskip

\item{[\MizGlobal]}
 N. Mizoguchi,
Boundedness of global solutions for a supercritical semilinear heat equation and its application,
{\sln Indiana Univ. Math. J.} 54 (2005), 1047-1059.
\smallskip

\item{[\NST]}
W.-M. Ni, P. Sacks and J. Tavantzis,
On the asymptotic behavior of solutions of certain quasilinear parabolic equations,
{\sln J. Differential Equations} 54 (1984), 97-120.
\smallskip

\item{[\OS]}
C. Ortner and E. S\"uli,
A note on linear elliptic systems on $\R^d$,
Preprint ArXiV 1202.3970.
\smallskip

\item{[\Pac]}
F. Pacard, Partial regularity for weak solutions of a nonlinear elliptic equation,
{\sln Manus\-cripta Math.}
79 (1993), 161-172.
\smallskip

\item{[\PQS]}
P. Pol\'a\v cik, P. Quittner and Ph. Souplet,
Singularity and decay estimates in superlinear problems
via Liouville-type theorems.  Part II: parabolic equations,
{\sln Indiana Univ. Math.~J.} 56 (2007), 879-908.
\smallskip
\item{[\PY]}
 P. Pol\'a\v cik and E. Yanagida,
On bounded and unbounded global solutions of a supercritical semilinear heat equation,
{\sln Math. Ann.} 327 (2003), 745-771.

\smallskip
\item{[\PYb]}
 P. Pol\'a\v cik and E. Yanagida,
Global unbounded solutions of the Fujita equation in the intermediate range,
{\sln Math. Ann.} 360 (2014), 255-266.
\smallskip

\item{[\Qui]}
P. Quittner,
A~priori bounds for global solutions of a semilinear parabolic problem,
{\sln Acta Math. Univ. Comenian. (N.S.)} 68 (1999), 195-203.

\smallskip
\item{[\QuiB]}
P. Quittner,
Threshold and strong threshold solutions of a semilinear parabolic equation,
Preprint (2016).
\smallskip

\item{[\QSb]} P. Quittner, Ph. Souplet,
Superlinear parabolic problems. Blow-up, global existence and steady states,
Birkhauser Advanced Texts, 2007, 584 p.+xi. ISBN: 978-3-7643-8441-8.
\smallskip

\item{[\SnTW]}
S. Snoussi, S. Tayachi and F.B. Weissler,
Asymptotically self-similar global solutions of a general semilinear heat equation,
{\sln Math. Ann.} 321 (2001), 131-155.
\smallskip

\item{[\SouCRAS]}
Ph.~Souplet,
Sur l'asymptotique des solutions globales pour une \'equation de la chaleur semi-lin\'eaire dans des domaines non born\'es, 
{\sln C. R. Acad. Sci. Paris S\'er. I Math.} 323 (1996), 877-882.
\smallskip

\item{[\SouCPDE]}
Ph.~Souplet,
Geometry of unbounded domains, Poincar\'e inequalities and stability 
in semilinear parabolic equations,
{\sln Comm. Partial Differential Equations}
24 (1999), 951-973.
\smallskip

\item{[\SouJFA]}
Ph.~Souplet,
Decay of heat semigroups in $L^\infty$ and applications 
to nonlinear parabolic problems in unbounded domains,
{\sln J. Funct. Anal.}
173 (2000), 343-360.
\smallskip

\item{[\Tay]}
 M.E. Taylor, 
 Analysis on Morrey spaces and applications to Navier-Stokes and other evolution equations,
 {\sln Comm. Partial Differential Equations} 17 (1992), 1407-1456.

\smallskip

\item{[\Tri]}
H. Triebel, Interpolation theory, function spaces, differential operators, 
North-Holland Mathematical Library 18, North-Holland, Amsterdam-New York, 1978. 

\smallskip
\item{[\YZ]}
M. Yamazaki and X. Zhou,
Semilinear heat equations with distributions in Morrey spaces as initial data,
{\sln Hokkaido Math. J.} 30 (2001), 537-571.
\smallskip

\item{[\Wei]}
F.B. Weissler,
Existence and non-existence of global solutions for a semilinear heat equation,
{\sln Israel J. Math.} 38 (1981), 29-40.

\bye